\documentclass[12pt]{iopart}
\usepackage{graphicx}
\usepackage{amssymb, amsmath}

\begin{document}

\title[$V$-variable fractals and superfractals]{$V$-variable fractals and superfractals }

\author{Michael Barnsley\dag, John E Hutchinson\dag\footnote[3]{To
whom correspondence should be addressed.} and 
\"{O}rjan Stenflo\ddag}

\address{\dag\ Department of Mathematics, Australian National University, Canberra, ACT. 0200, Australia}
\address{\ddag\ Department of Mathematics, Stockholm University, SE-10691, Stockholm, Sweden}

\begin{abstract}
Deterministic and random fractals, within the framework of Iterated Function Systems, have been used to model and study a wide range of phenomena across many areas of science and technology.  However, for many applications deterministic fractals are locally too similar near distinct points while standard random fractals have too little local correlation.  Random fractals are also slow and difficult to compute.  These two major problems restricting further applications are solved here by the introduction of  $V$-variable fractals and superfractals.  
\end{abstract}


\maketitle

\section{Introduction}
	
We introduce the class of  $V$-variable fractals.  This solves two major problems previously restricting further applications of Iterated Function System (IFS) generated fractals, the absence of fine control on local variability and the absence of a fast algorithm for computing and accurately sampling standard random fractals. 

The integer parameter  $V$  controls the number of distinct shapes or forms at each level of magnification (Figures 4, 6). The case  $V=1$ includes standard deterministic fractals (Figures 1, 2, also the fern and lettuce in Figure 4) generated by a single IFS and homogeneous random fractals  (c.f.\ Hambly 1992, 2000 and Stenflo 2001).  Large  $V$ allows rapid approximations to standard random fractals in a quantifiable manner, and the approximation is to not one, but to a potentially infinite sequence of correctly distributed examples. The construction of  $V$-variable fractals is by means of a fast Markov Chain Monte Carlo type algorithm.   In particular, one can now approximate stndard random fractals, together with their associated probability distribution, by means of a fast forward algorithm.     

A surprising but important fact is that each family of  $V$-variable fractals, together with its naturally associated probability distribution, forms a single ÒsuperÓfractal generated by a single ÒsuperÓIFS, operating not on points in the plane (for example) as for a standard IFS, but on  $V$-tuples of images.  Furthermore, dimensions of  $V$-variable fractals are computable using products of random matrices and ideas from statistical mechanics. We implement a Monte Carlo method for this purpose. 

Since the mathematically natural notion of a $V$-variable fractal solves two major problems previously restricting wider applications, we anticipate that  $V$-variable fractals should lead to significant developments and applications of fractal models.  For example, IFSs provide models for certain plants, leaves, and ferns, by virtue of the self-similarity which often occurs in branching structures in nature.  But nature also exhibits randomness and variation from one level to the next --- no two ferns are exactly alike, and the branching fronds become leaves at a smaller scale.  $V$-variable fractals allow for such randomness and variability across scales, while at the same time admitting a continuous dependence on parameters which facilitates geometrical modelling. These factors allow us to make the hybrid ÒbiologicalÓ models in Figures 4, 9.  Because of underlying special code trees (Barnsley, Hutchinson and Stenflo 2003a) which provide the fundamental information-theoretic basis of  $V$-variable fractals, we speculate that when a  $V$-variable geometrical fractal model is found that has a good match to the geometry of a given plant, then there is a specific relationship between these code trees and the information stored in the genes of the plant.

In this paper we describe the algorithm for generating $V$-variable fractals, explain the notion of a superfractal, and show how to compute the dimension of  $V$-variable fractals in case a uniform open set condition applies.  In order to make the ideas clearer, and hopefully facilitate the application of these notions to non-mathematical areas, we have avoided technicalities and  illustrated the ideas by means  of  a number of model examples, but the general constuction and results should be clear.   

In Barnsley, Hutchinson and Stenflo (2003a) we survey the classical properties of IFSs, develop the underlying theory of V-variable code trees and establish general  existence and other properties for $ V $-variable fractals and superfractals.  In Barnsley, Hutchinson and Stenflo (2003b) we prove the dimension results for $V$-variable fractals, for which we here give an informal justification.

\section{Iterated Function Systems}
	By way of background we first recall the concept of an IFS via the canonical example of the Sierpinski triangle  $S$ (approximated in the bottom right panel of Figure 1), which has been studied extensively (see Falconer 1990, 1997 and Hambly 1992, 2000) both mathematically and as a model for diffusion processes through disordered and highly porous material.   The set $S$ has three components, each of which is a scaled image of itself; each of these components has three sub-components, giving nine scaled images of   at the next scale, and so on ad infinitum.  A simple observation is that if  $f_1, f_2, f_3 $  are contractions of space by the factor $ \frac{1}{2} $  with fixed points given by the three vertices   
$A_1, A_2, A_3$ respectively of  $S$, then the three major Òsub-componentsÓ of  $S$ are  $f_1(S), f_2(S), f_3(S) $  respectively and  
\begin{equation}\label{ffp}
 S = f_1(S) \cup f_2(S) \cup f_3(S). 
 \end{equation}
 
 \begin{center}
\includegraphics[scale=1]{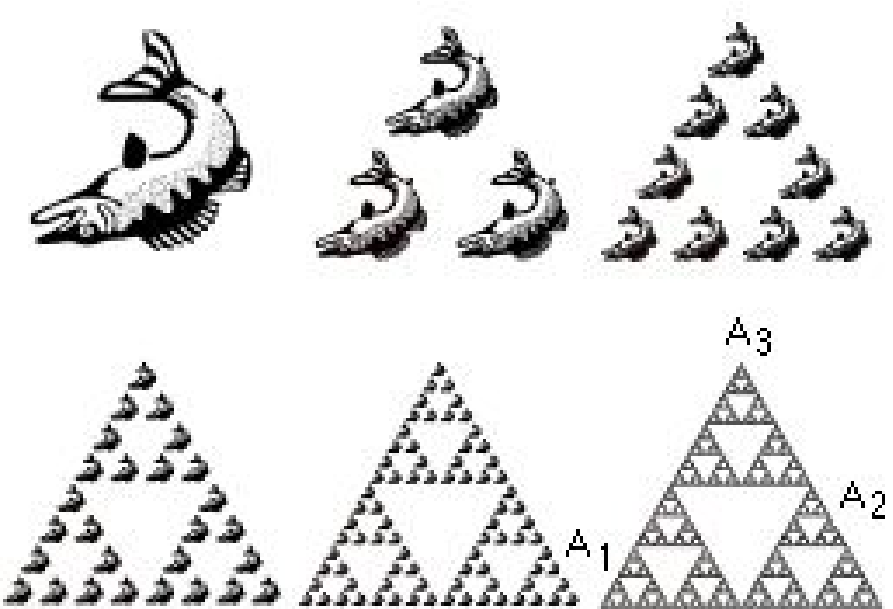}\\
\end{center}
\begin{quote}
\textsf{\textbf{Figure 1.}  Convergent or Backward Process.  Beginning from any set (fish) $T_0$,  iterates   $T_1=F(T_0)$, $T_2=F(T_1)$,\dots converge to the Sierpinski Triangle $S$.  Shown are iterates $T_0,T_1,T_2,T_3,T_4,T_8$.}
\end{quote}

The collection of maps  $F=(f_1,f_2,f_3) $  is called an \emph{Iterated Function System} (or \emph{IFS}).  For any set  $T$  one similarly defines  $F(T) = f_1(T) \cup f_2(T) \cup f_3(T)$.  It is not difficult to show that  $S$ is the unique compact set satisfying
\eqref{ffp} .  Furthermore, beginning from any compact set  $T_0$, and for  $k\geq 1 $ recursively defining $T_k = F(T_{k-1})$, it follows that   $T_k$ converges to $S$ in the Hausdorff metric  as  $k\to \infty$, independently of the initial set  $T_0$ (Figure 1).  For this reason,  $S$ is called the \emph{fractal set attractor} of the IFS $F$ and this approximation method is called the \emph{convergent} or \emph{backward process}, c.f.\ Hutchinson (1981).

An alternative approach to generating  $S$  is by a \emph{chaotic} or \emph{forward process} (Barnsley and Demko 1985), sometimes called the \emph{chaos game} (Figure 2).  Begin from any point   $x_0 $ in the plane and recursively define $x_k=\widehat{f}_k(x_{k-1})$, where each  
$\widehat{f}_k$ is chosen independently and with equal probability from 
$(f_1,f_2,f_3)$.  With probability one the sequence of points  $(x_k)_{k\geq 0}$ approaches and moves ÒergodicallyÓ around, and increasingly closer to, the attractor  $S$.  For this reason  $F$ is called an \emph{iterated} function system.  If instead the   $\widehat{f}_k$ are selected from  $(f_1,f_2,f_3)$ with probabilities $(p_1,p_2,p_3)$  respectively, where each $p_i>0 $   and $p_1+p_2+p_3=1$, then the same set  $S$ is determined by the sequence $(x_k)_{k\geq 0}$, but now the points accumulate unevenly, and the resulting \emph{measure attractor} can be thought of as a greyscale image on  $S$, or probability distribution on  $S$, or more precisely as a measure.  In this case  $(f_1,f_2,f_3;p_1,p_2,p_3)$ is called an \emph{IFS with weights}.

\begin{center}
\includegraphics[scale=1]{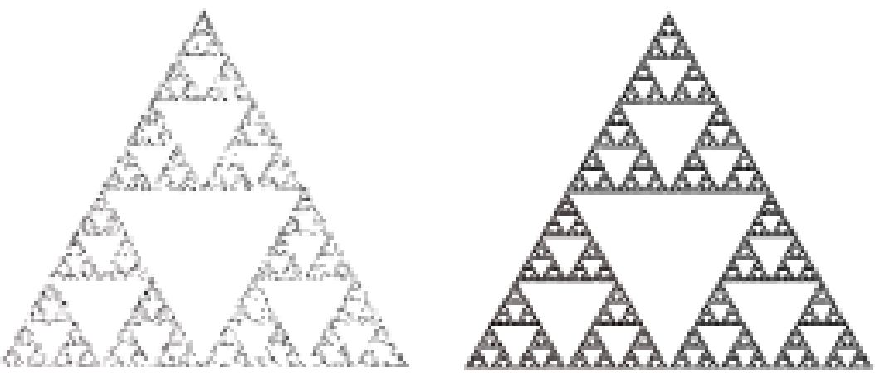}
\end{center} 
\begin{quote}
\textsf{\textbf{Figure  2.}  Chaotic or Forward Process.  Beginning from any initial point and randomly and independently applying $f_1$, $f_2$ or $f_3$   produces the Sierpinski Triangle   as attractor with probability one. Shown are the first 10,000 and 100,000 points respectively.}
\end{quote}

These ideas and results naturally extend to general families of contraction maps and probabilities.  Even with a few affine or projective transformations, one can construct natural looking images (see the initial fern and lettuce in Figure 4).  

IFSs have been extended to study the notion of random fractals. See Falconer (1986), Graf (1987) and Mauldin and Williams (1986); also Hutchinson and R\"{u}schendorf (1998, 2000) where the  idea of an IFS operating directly in the underlying probability spaces is used.

\section{Construction of $ V $-variable fractals}

We now proceed to the construction of $V$-variable fractals.  This can be understood in the model situation of two IFSs  $F$ and $G$, with $F$    as before and with  $G$ having the same fixed points  $A_1,A_2,A_3$  but with contraction ratios $ \frac{1}{3} $  instead of $ \frac{1}{2}$. We emphasise that the following construction is in no sense ad hoc, but is the natural chaotic or forward process for a \emph{superfractal} whose members are  $V$-variable fractals as we see later.  

\begin{center}
\includegraphics[scale=1]{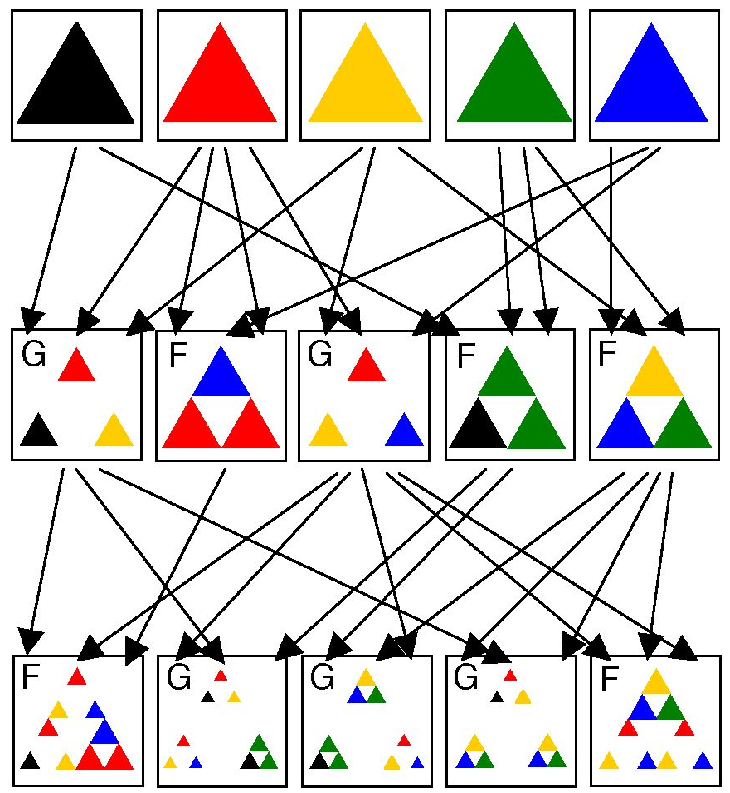}
\end{center}
\begin{quote}
\textsf{\textbf{Figure  3.}   Forward algorithm for generating families of  $V$-variable fractals.  Shown are levels 1 (top), 2 and 3 in the construction of a potentially infinite sequence of 5-tuples of 5-variable Sierpinski Triangles.  For each buffer from level 2 onwards, F or G indicates which IFS was used, and the input arrows indicate the buffers to which this IFS was applied.}
\end{quote}

One begins with arbitrary sets, one in each of  $V$ input buffers at level 1 (Figure 3 where $V=5$, colour coded to indicate the maps involved, and Fig. 6 where $V=2$).  A set in the first of  $V$ output buffers (level 2, Figure 3) is constructed as follows: choose an IFS  $F$ or $G$  with the preassigned probabilities  $P^F$ or  $P^G$ respectively; then apply the chosen IFS to the content of three buffers chosen randomly and independently with uniform probability from the input buffers at level 1, allowing the possibility that the same buffer is chosen more than once  (thus one is performing uniform sampling with replacement).   The resulting set is placed in the first buffer at level 2.   The content of each of the remaining buffers at level 2 is constructed similarly and independently of the others at level 2.  These output buffers then become the input buffers for the next step and the process is repeated, obtaining the bottom row in Figure 3, and so on. 

The construction produces an arbitrarily long sequence of  $V$-tuples of \emph{approximate} $V$-variable fractals associated to the pair of IFSs $(F,G)$  and the probabilities  $(P^F,P^G)$.  The degree of approximation is soon, and thereafter remains, within screen resolution or machine tolerance (Figure 6).  The empirically obtained distribution of $V$-variable fractals over any infinite sequence of runs is the same with probability one and is the natural distribution as we explain later.  The generalisation to the case of a family of IFSs   $(F^1,\dots,F^N)$ with associated probabilities $(P^1,\dots, P^N)$  is straightforward; also sets can be replaced by greyscale images or more generally by coloured sets built up from primary colours if IFSs with weights are used (Figure 9). 

\begin{center}
\includegraphics[scale=1]{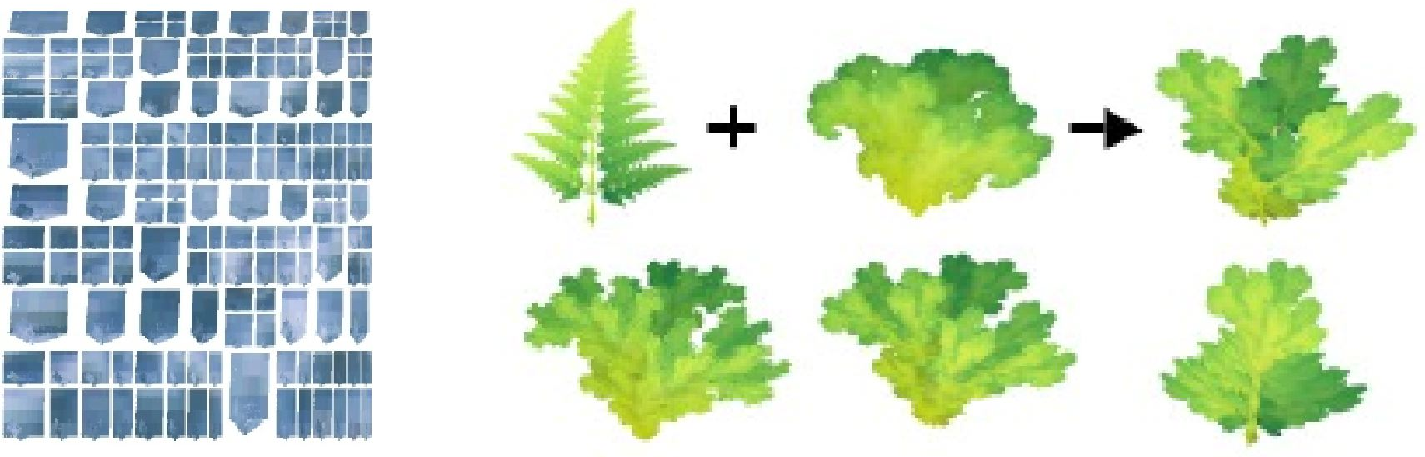}
\end{center}
\begin{quote}
\textsf{\textbf{Figure  4.}   2-Variable Fractals.  In the first ``Squares and Shields'' example, the number of components at each level of magnification is 1,4,16,64, \dots, but at each level there are at most two distinct shapes/forms up to affine transformations. The actual shapes depend on the level.  In the second example of Òfractal breedingÓ, fern and lettuce fractal parents are shown with four possible 2-variable offspring. The two IFSs used are those generating the fern and the lettuce. The associated superfractal is the family of all possible offspring together with the naturally associated probability distribution.}
\end{quote}

The \emph{$V$-variability} can be understood as follows.  In Figure 3 the set in each buffer from level 2 onwards is composed of three component parts, each obtained from one of the  $V=5$ buffers at the previous level; at level 3 onwards sets are composed of 9 smaller component parts each obtained from one of the  $V$ buffers two levels back; at level 4 onwards sets are composed of 27 smaller component parts each obtained from one of the $V$  buffers three levels back; etc.  Thus at each level of magnification there are at most  $V$ distinct component parts up to rescaling.   In general, for a  $V$-variable fractal, although the number of components grows exponentially with the level of magnification, the number of distinct shapes or forms is at most  $V$ up to a suitable class of transformations (e.g.\ rescalings, affine or projective maps) determined by the component maps in each of the IFSs  $(F^1,\dots,F^N)$ (Figure 4 first panel, Figure 6).  If parts of the  $V$-variable fractals overlap each other, the $V$-variability is not so obvious, and fractal measures (greyscale or colour images) rather than fractal sets (black and white images) are then more natural to consider.

\section{Superfractals}
The limit probability distribution on the infinite family of  $V$-variable fractals obtained by the previous construction is independent of the experimental run with probability one as we explain later.  As noted previously, an initially surprising but basic fact is that this family of  $V$-variable fractals and its probability distribution is a fractal in its own right, called a \emph{superfractal}, and the construction process for the generated collection of  $V$-variable fractals turns out to be the forward or chaotic process for this superfractal, see also Barnsley, Hutchinson and Stenflo (2003a).  For large  $V$,  $V$-variable fractals approximate standard random fractals and their naturally associated probability distributions in a quantifiable manner, providing another justification for the canonical nature of the construction, see Barnsley, Hutchinson and Stenflo (2003a, b).  Although there was previously no useful forward algorithm for standard random fractals such as Brownian sheets, one can now rapidly generate correctly distributed families of such fractals to any specified degree of approximation by using the previously described fast forward (Monte Carlo Markov Chain) algorithm for large  $V$.  The number of required operations typically grows linearly in  $V$ as only sparse matrix type operations are required.  
  
The superfractal idea can be understood as follows.  The process of passing from the  $V$-tuple of sets at one level of construction to the  $V$-tuple at the next (Figures 3, 6) is given by a random function $\mathcal{F}^a : \mathbf{H}^V \to \mathbf{H}^V$, where   $\mathbf{H}^V $ is the set of all $V$-tuples of compact subsets of the plane $\mathbb{R}^2 $  or of some other compact metric space as appropriate, and where   $a$ belongs to some index set  $\mathcal{A}$.  All information necessary to describe the chosen $\mathcal{F}^a$ at each level  in Figure 3 is given by the chosen IFSs for each buffer at that level (namely $G,F,G,F,F$   across level 2) and by the arrows pointing to each buffer at that level (indicating which three buffers at the previous level are used for each application of  $F$ or of $G$),  see also the captions in Figure 6.  Each $\mathcal{F}^a$  has a certain probability  $\mathcal{P}^a$ of being chosen, this probability is induced in the natural manner from the probabilities  $P^F$  and  $P^G$  of selecting $F$  or  $G$.  The  $\mathcal{F}^a$ are contraction maps on  $\mathbf{H}^V $ in the Hausdorff metric, with contraction ratio equal $\frac{1}{2} $, and in general the contraction ratio of $\mathcal{F}^a$ equals the maximum of the contraction ratios of  the individual maps in the IFSs being used.  In particular,  the superIFS  $(\mathcal{F}^a, \mathcal{P}^a, a\in \mathcal{A})$  is an IFS operating not on points in  $\mathbb{R}^2 $ as for a standard IFS but on $V$-tuples of sets in $\mathbf{H}^V $.  From IFS theory applied in this setting, there is a unique superfractal set and superfractal measure which with probability one is effectively given by the collection of  $V$-tuples of  $V$-variable fractals together with the experimentally obtained probability distribution arising from the previous construction. See Barnsley, Hutchinson and Stenflo (2003a) for detailed proofs.

\section{Examples  of $ 2 $-variable fractals}

We begin with 2 IFSs $ U=(f_1,f_2) $ (``Up with a reflection'') and $ D=(g_1,g_2) $ (``Down''), where
\begin{align*}
f_1(x,y) &= \Big( \frac{x}{2} +\frac{3y}{8}- \frac{1}{16},
        \phantom{+}\frac{x}{2}-\frac{3y}{8}+\frac{9}{16}\Big),
&f_2(x,y)  & =        \Big( \frac{x}{2} -\frac{3y}{8}+ \frac{9}{16},
         -\frac{x}{2}-\frac{3y}{8}+\frac{17}{16}\Big),\\
g_1(x,y) &= \Big( \frac{x}{2} +\frac{3y}{8}- \frac{1}{16},
         -\frac{x}{2}+\frac{3y}{8}+\frac{7}{16}\Big),
&g_2(x,y)  & =        \Big( \frac{x}{2} -\frac{3y}{8}+ \frac{9}{16},
         \phantom{+}\frac{x}{2}+\frac{3y}{8}-\frac{1}{16}\Big).
    \end{align*}
The corresponding fractal attractors are shown in Figure 5.

\begin{center}
\includegraphics[scale=.2]{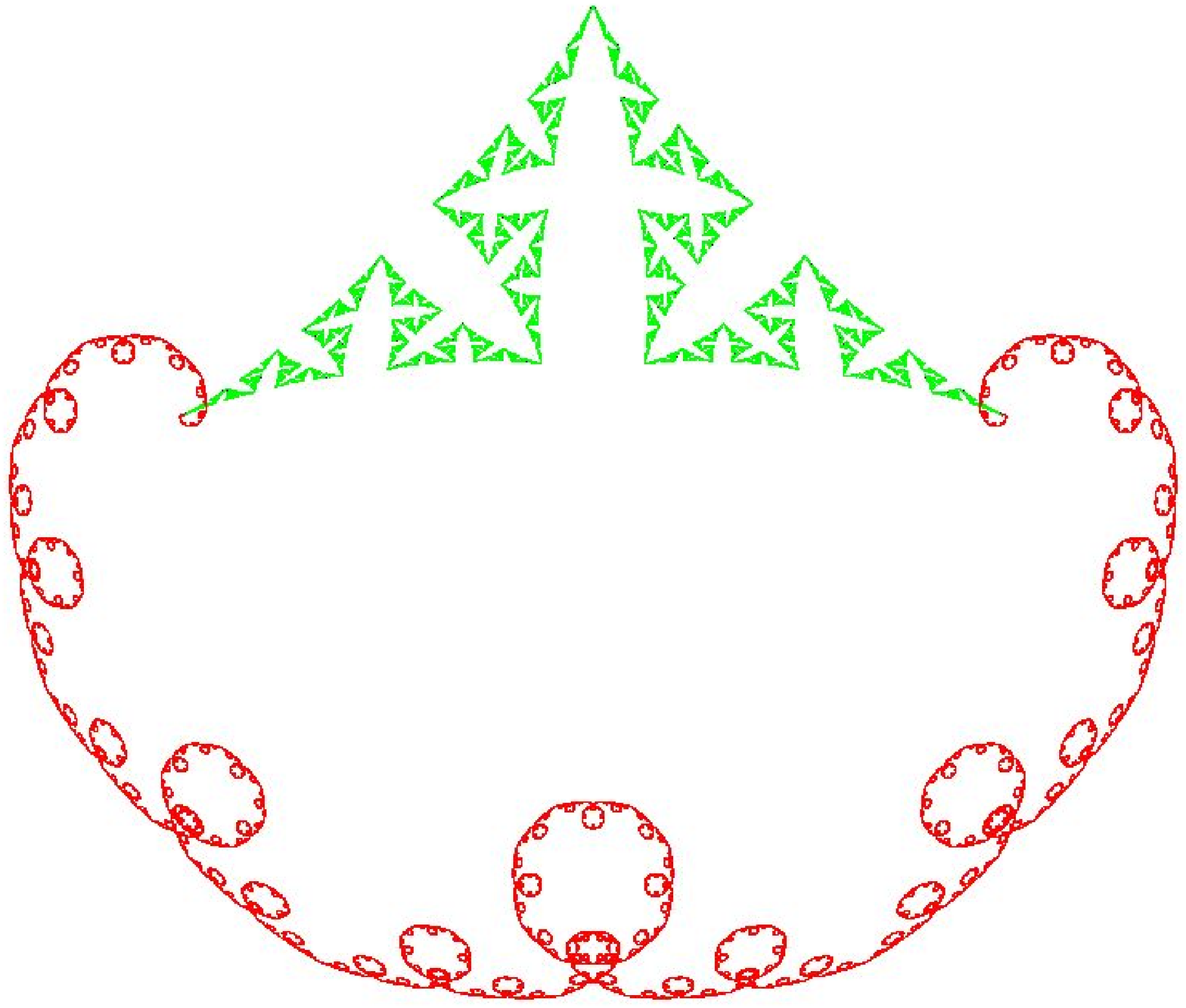} 

\textsf{\textbf{Figure  5.}  Up  (green) and  Down (red) attractors}
\end{center}

\medskip

Figure 6 shows the first 20 steps in the construction of a sequence of pairs of 2-variable fractals from the two IFSs $U $ and $ D $.  The initial pair of input figures can be arbitrarily chosen, here they are each the same and consist of four leaves.  

For the first step in the construction (producing the contents of the second pair of buffers) the IFS $ U =(f_1,f_2)$ was chosen, $ f_1$  was applied to  the previous left buffer $ L $ and $ f_2 $ was applied to the previous  right buffer  $ R $; the second buffer was obtained by applying $ U $ with $ f_1 $ and $ f_2 $ both  acting on the right buffer $ R $ at the previous step. Thus the first step in the construction can be described  by $ U(L,R) $ and $ U(R,R) $ respectively; see the caption below the second pair of screens in Figure 6.  The second step is given by $D(R,R)$ and $D(R,L)$, the third by $U(L,R)$ and $D(R,L)$, the fourth by $U(R,R)$ and $D(L,L)$, and so on from left to right and then down the page.

\bigskip

\begin{center}
\includegraphics[scale=.09]{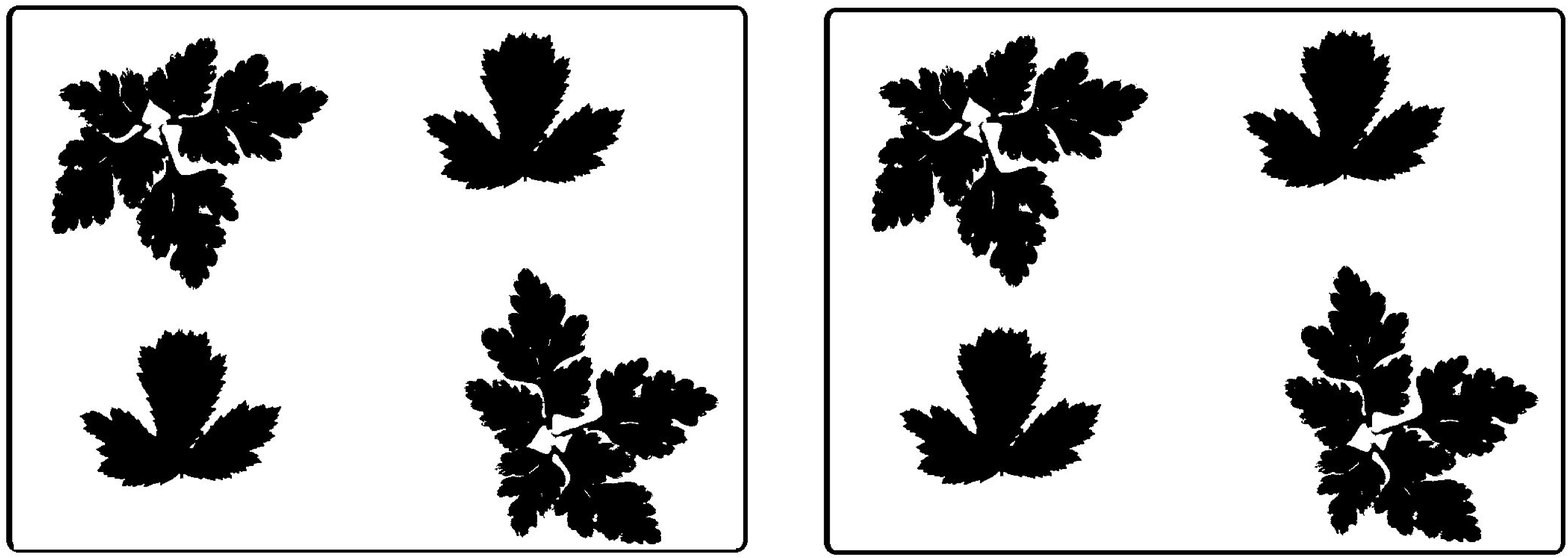}\\
\textsf{Initial Sets}
\end{center}

\begin{center}
\includegraphics[scale=.09]{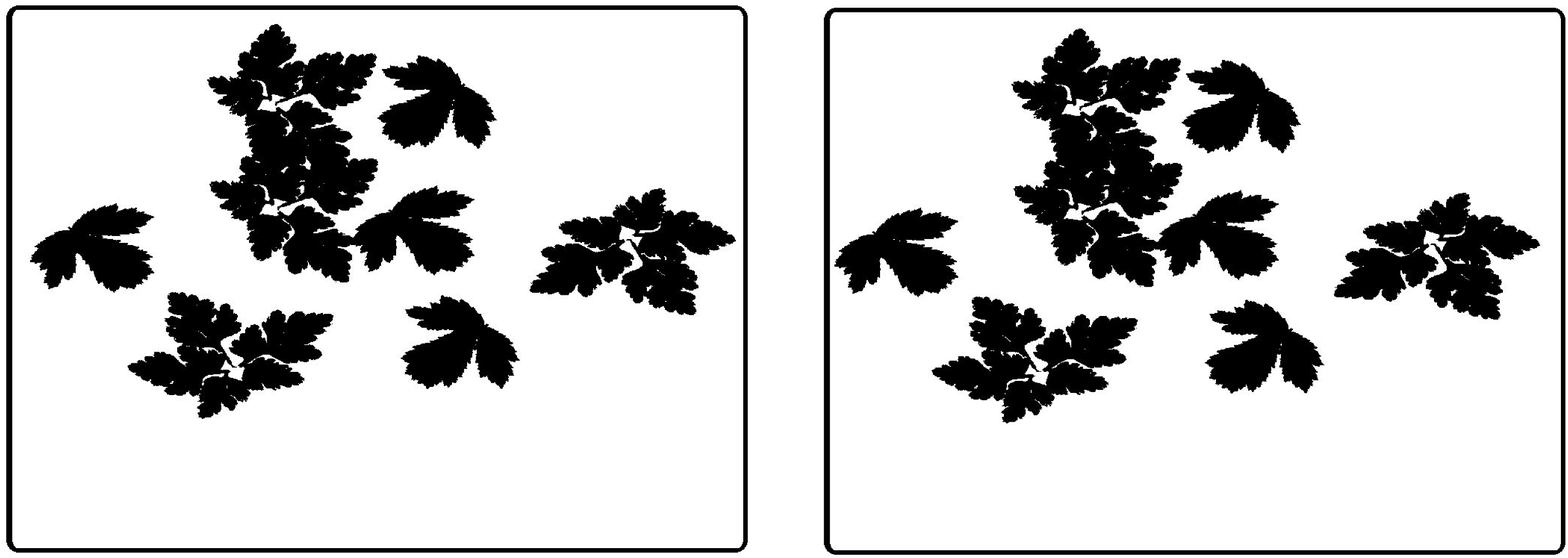} \qquad  
\includegraphics[scale=.09]{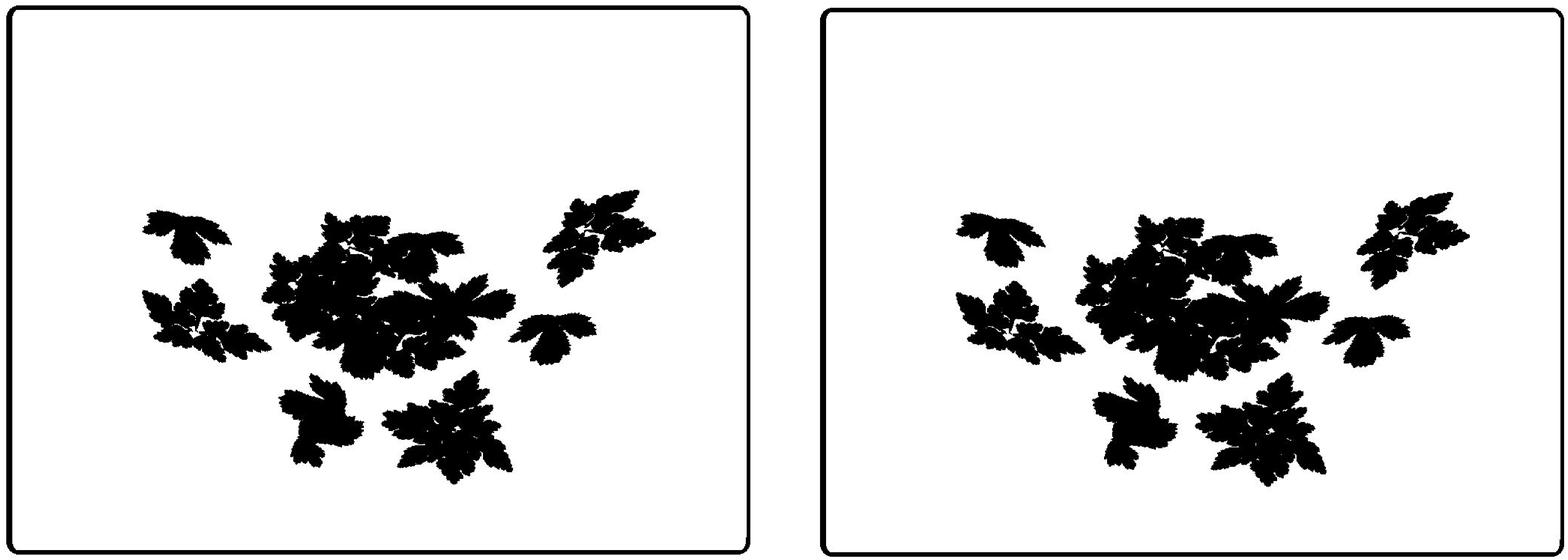}\\
$U(L,R) $  \hspace{1.9cm} $U(R,R) $ \hspace{2.6cm} $D(R,R) $   \hspace{1.9cm} $D(R,L)$
\end{center} 

\begin{center}
\includegraphics[scale=.09]{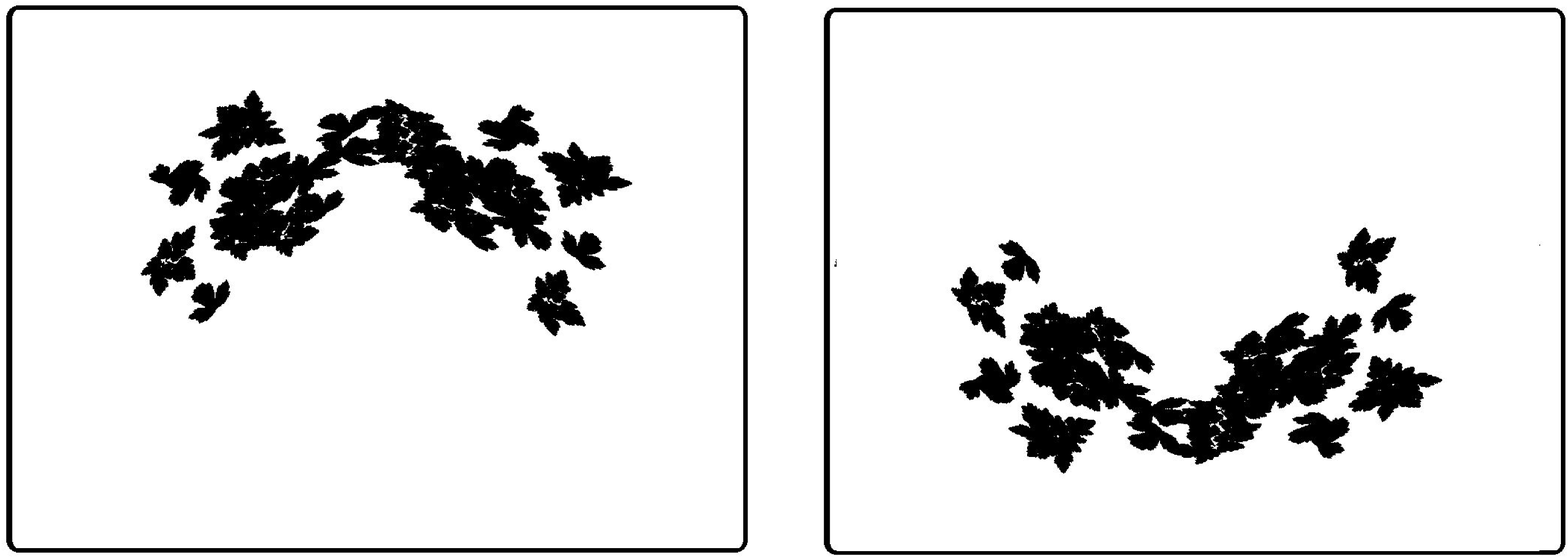}  \qquad  
\includegraphics[scale=.09]{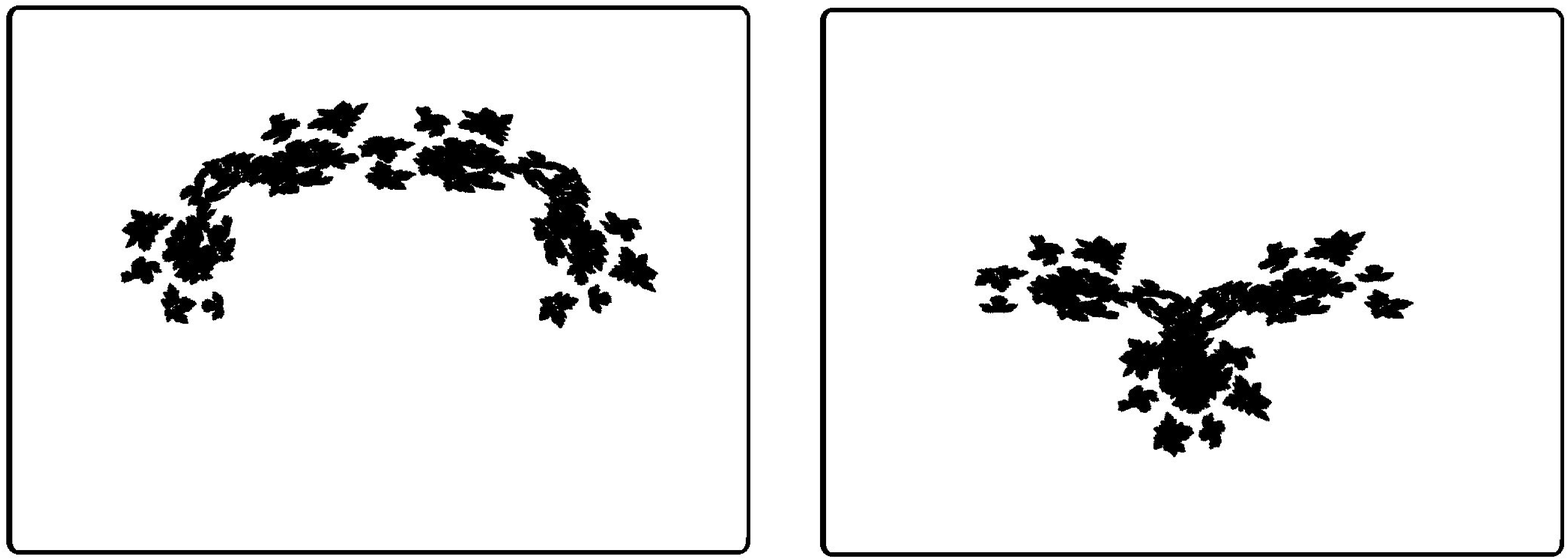}\\
$U(L,R) $  \hspace{1.9cm} $D(R,L) $ \hspace{2.6cm} $U(R,R) $   \hspace{1.9cm} $D(L,L)$
\end{center} 

\begin{center}
\includegraphics[scale=.09]{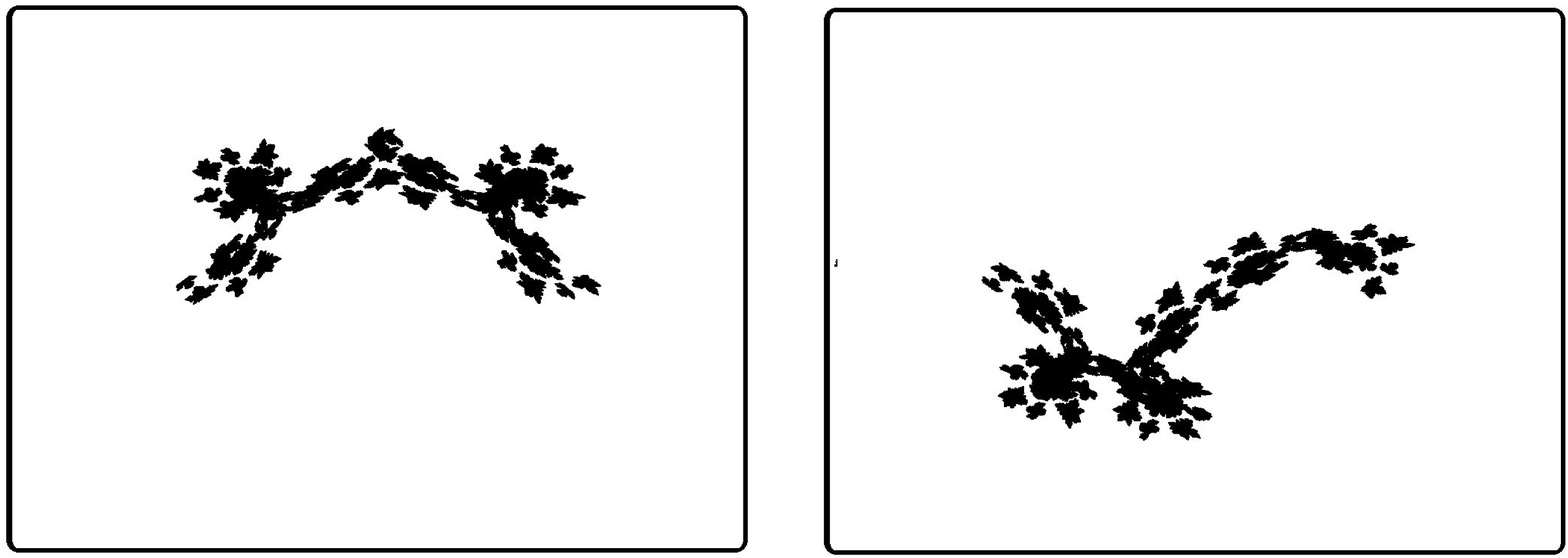}  \qquad  
\includegraphics[scale=.09]{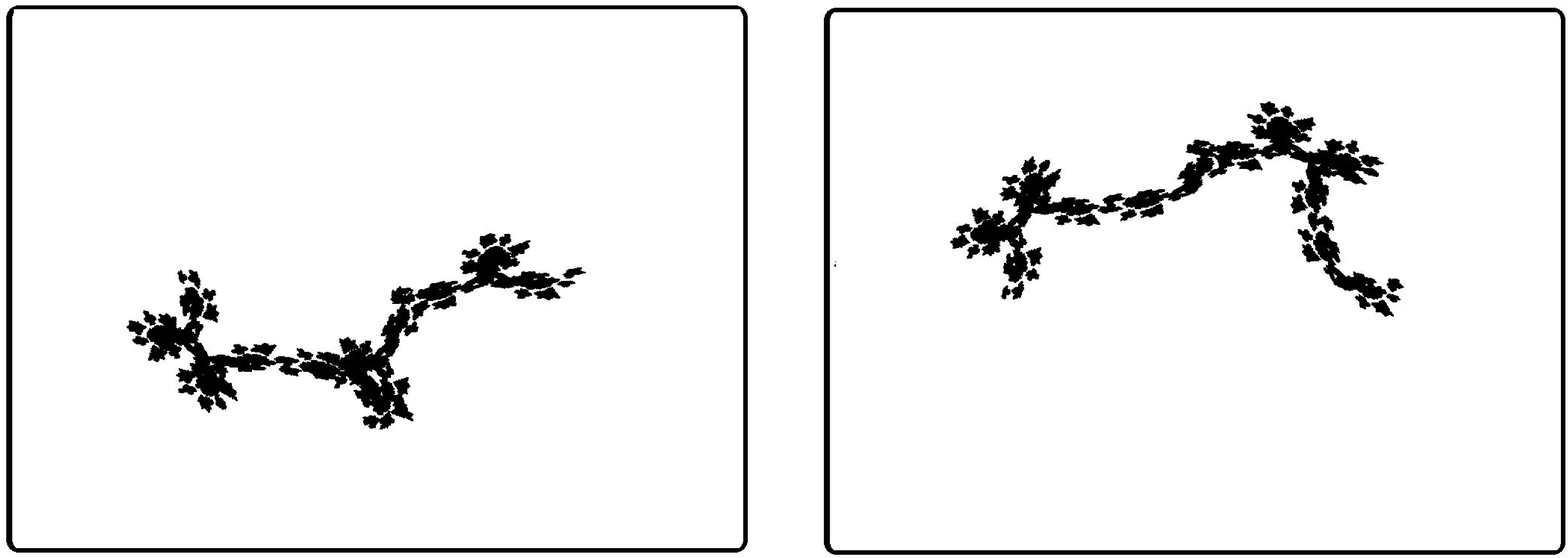} \\
$U(R,R) $  \hspace{1.9cm} $D(L,R) $ \hspace{2.6cm} $D(R,L) $   \hspace{1.9cm} $U(R,R)$
\end{center} 

\begin{center}
\includegraphics[scale=.09]{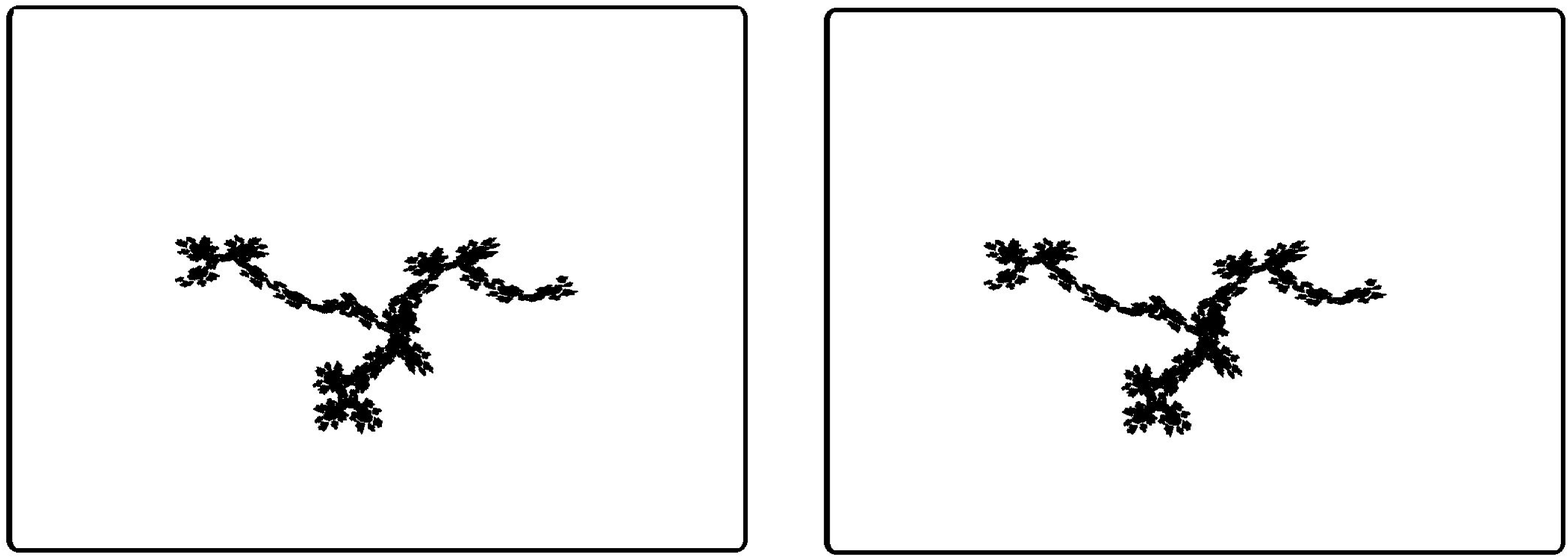} \qquad  
\includegraphics[scale=.09]{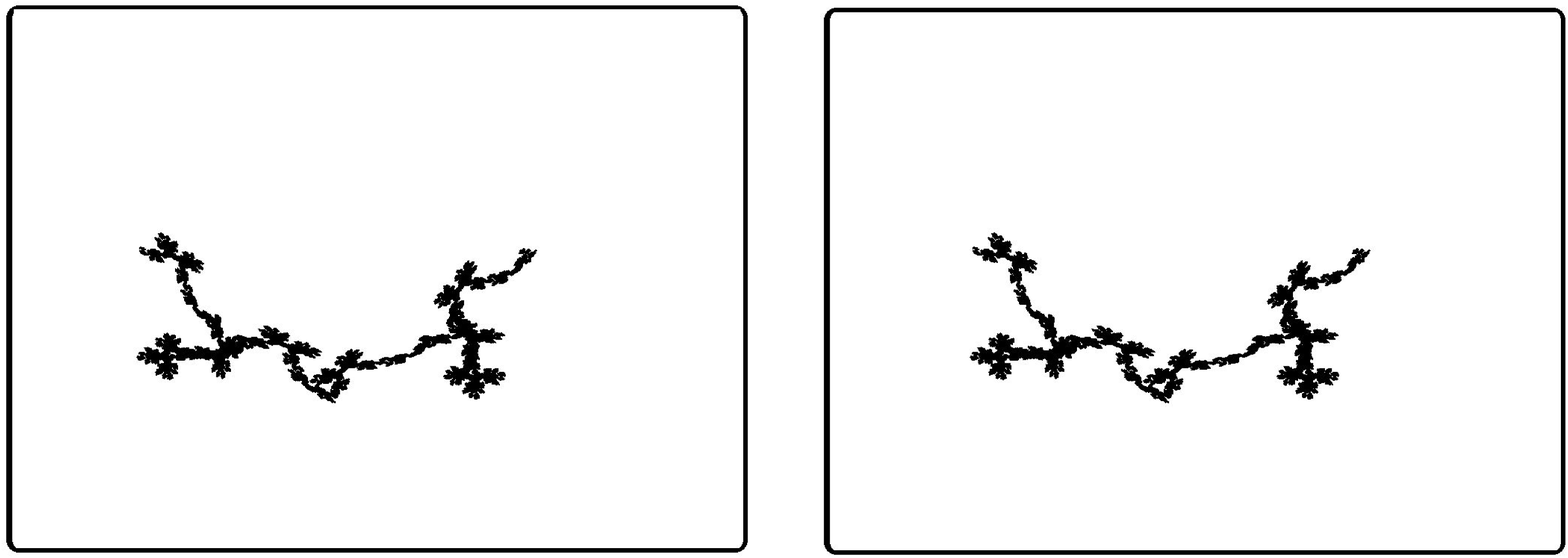}\\ 
$D(R,R) $  \hspace{1.9cm} $D(R,R) $ \hspace{2.6cm} $D(L,L) $   \hspace{1.9cm} $D(R,L)$
\end{center} 

\begin{center}
\includegraphics[scale=.09]{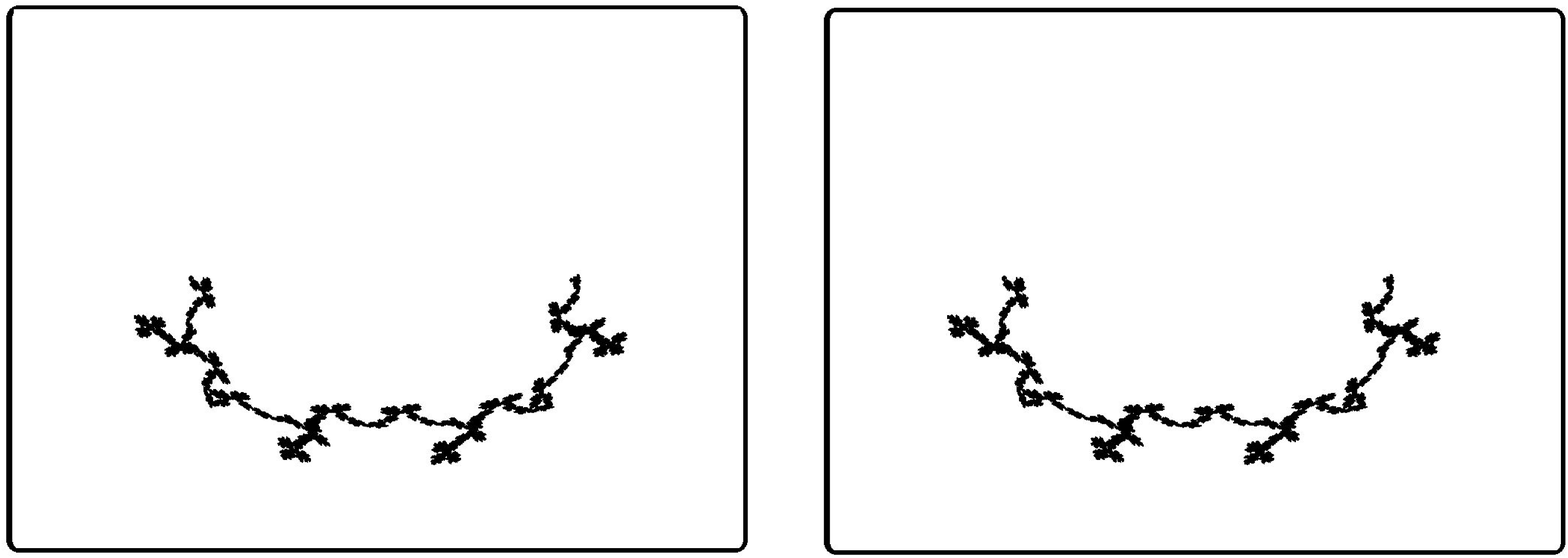}  \qquad  
\includegraphics[scale=.09]{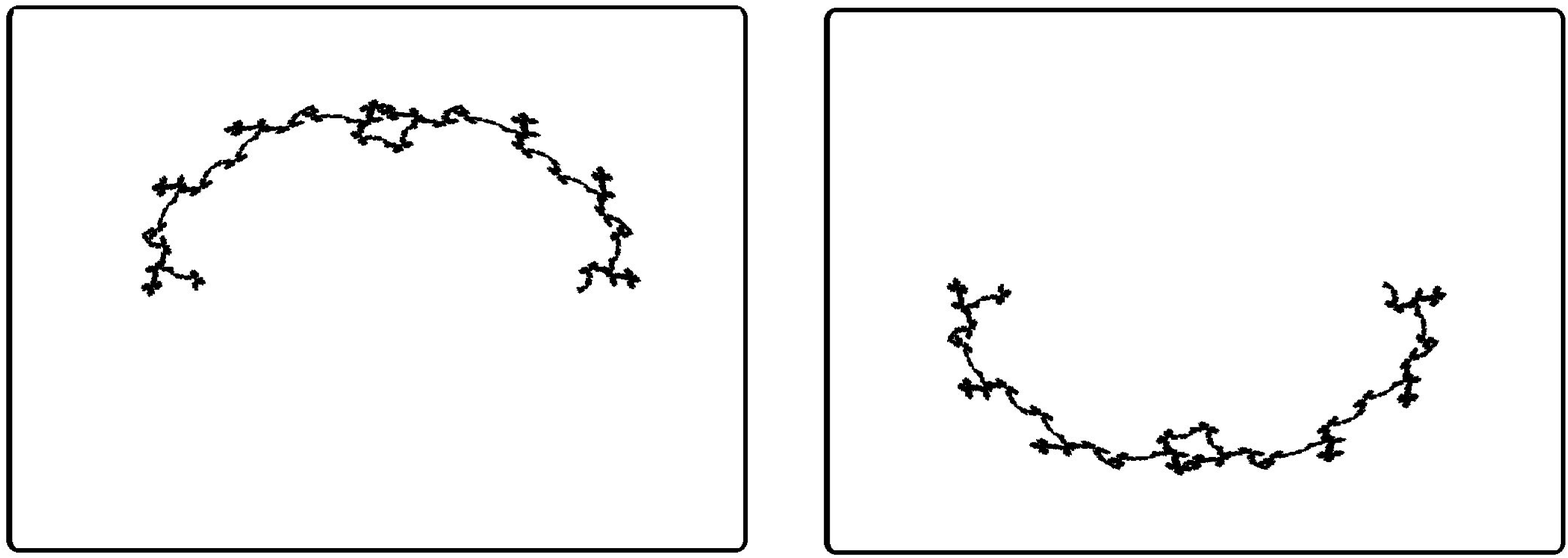}\\ 
$D(L,R) $  \hspace{1.9cm} $D(L,L) $ \hspace{2.6cm} $U(R,R) $   \hspace{1.9cm} $D(R,L)$
\end{center} 

\begin{center}
\includegraphics[scale=.09]{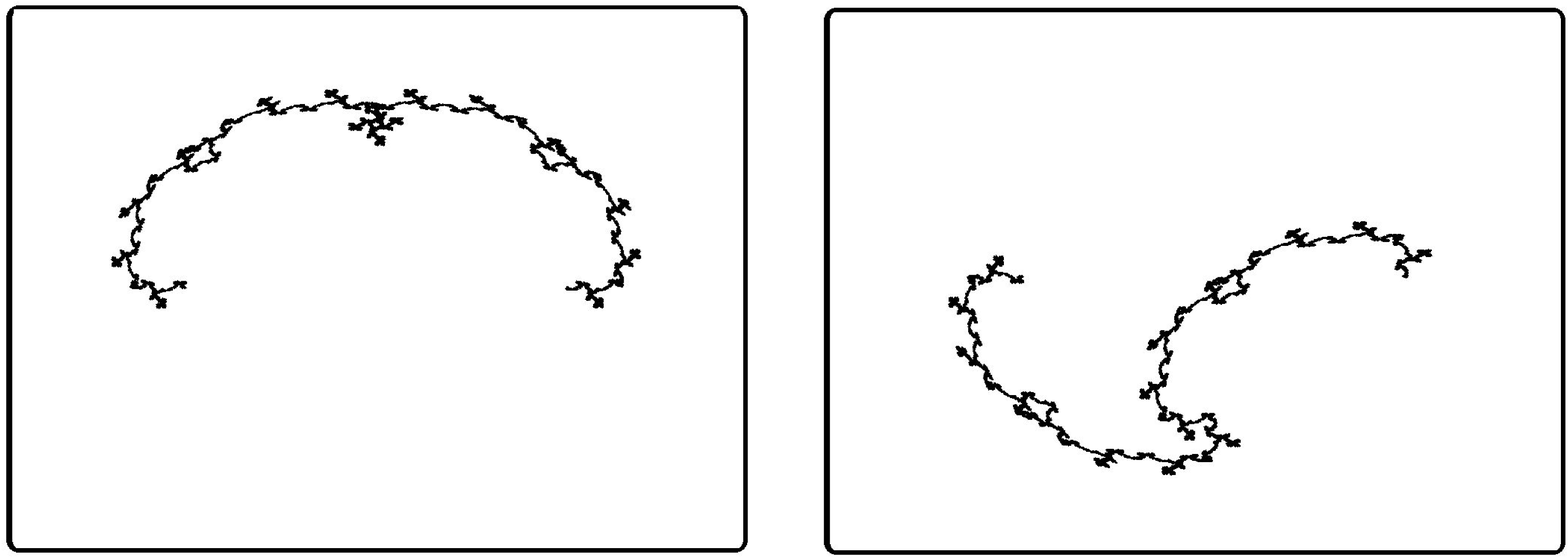} \qquad  
\includegraphics[scale=.09]{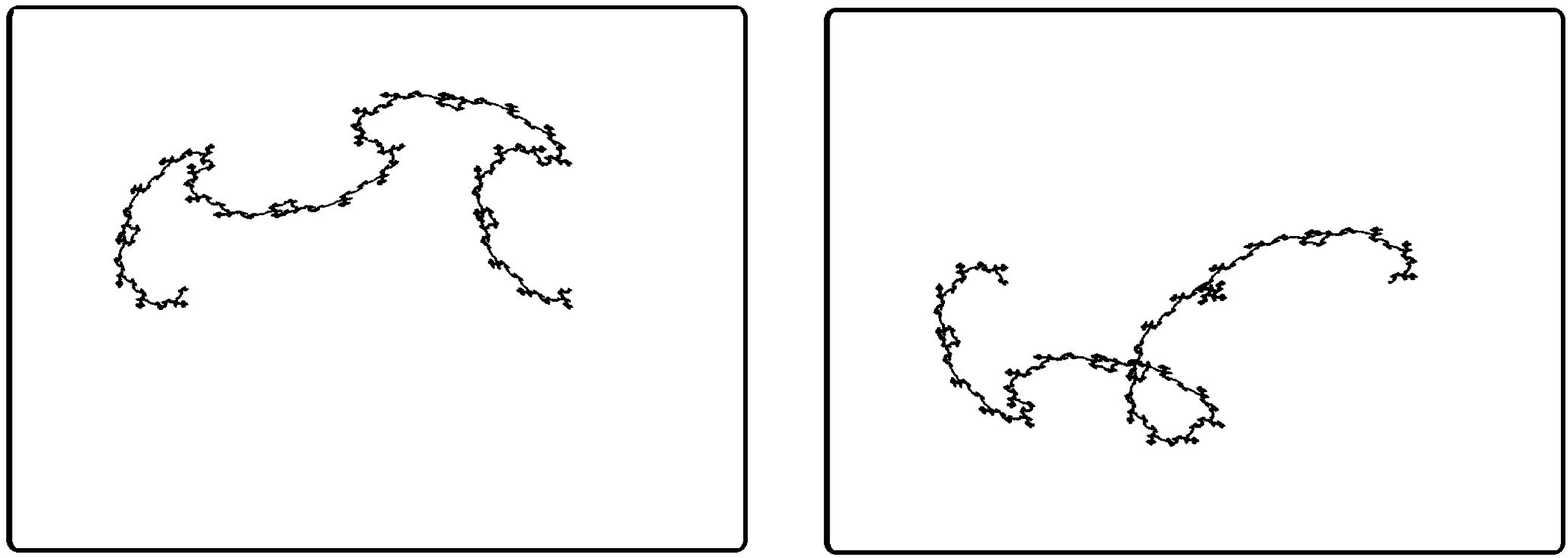}\\ 
$U(R,R) $  \hspace{1.9cm} $D(R,L) $ \hspace{2.6cm} $U(R,R) $   \hspace{1.9cm} $D(R,L)$
\end{center} 

\begin{center}
\includegraphics[scale=.09]{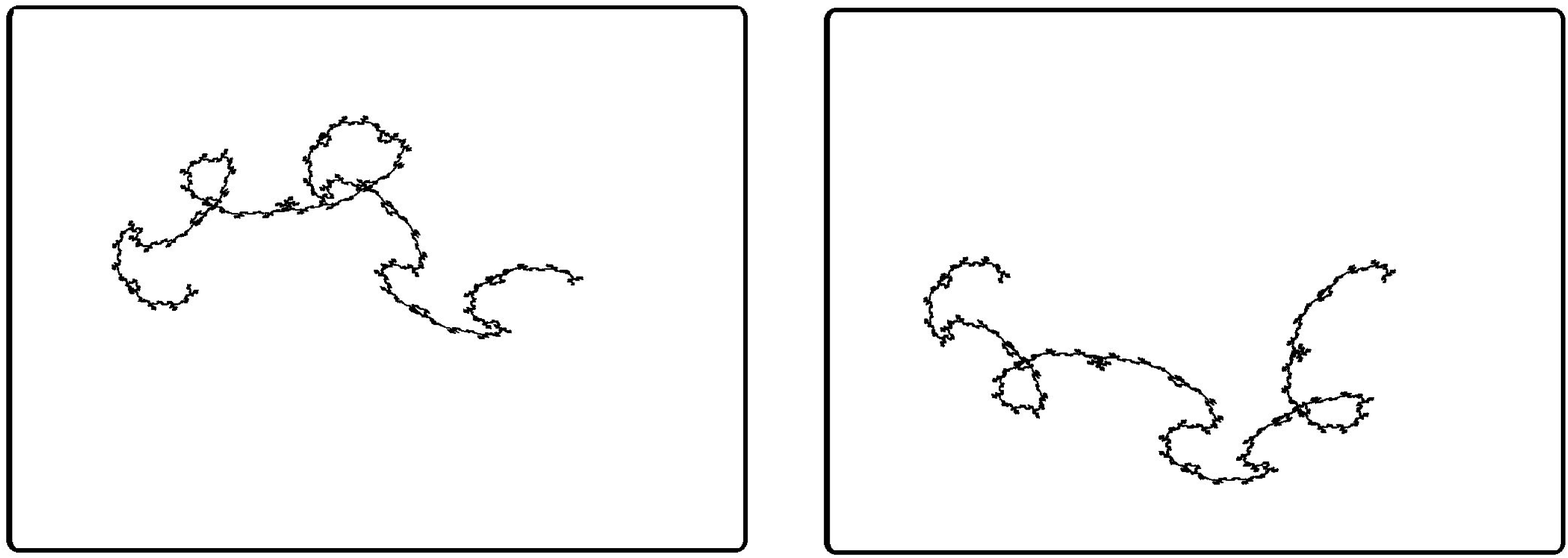}  \qquad  
\includegraphics[scale=.09]{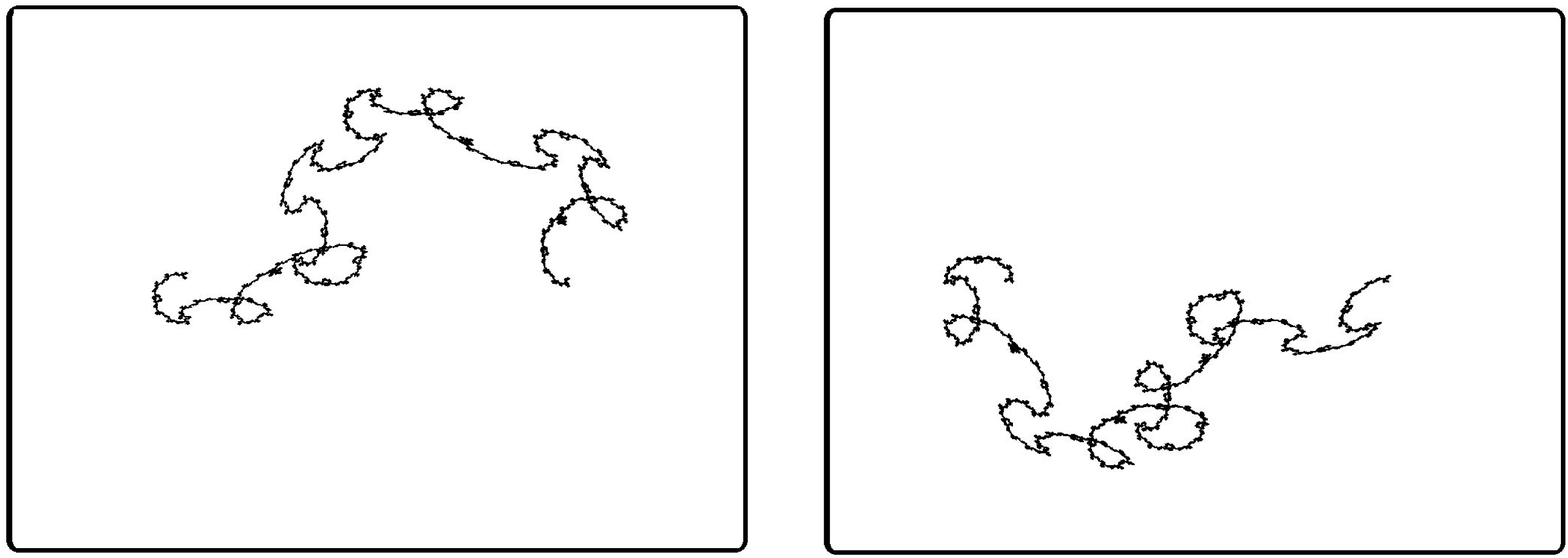}\\ 
$U(R,L) $  \hspace{1.9cm} $D(R,R) $ \hspace{2.6cm} $U(L,R) $   \hspace{1.9cm} $D(R,L)$
\end{center} 

\begin{center}
\includegraphics[scale=.09]{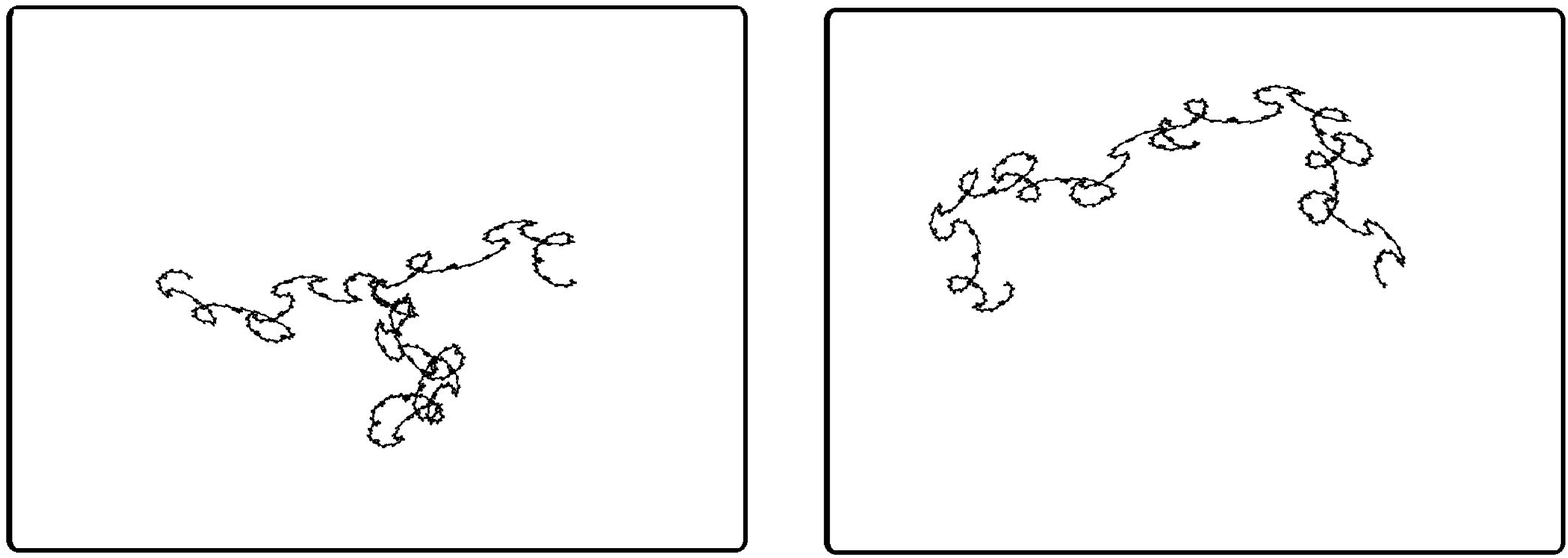}  \qquad  
\includegraphics[scale=.09]{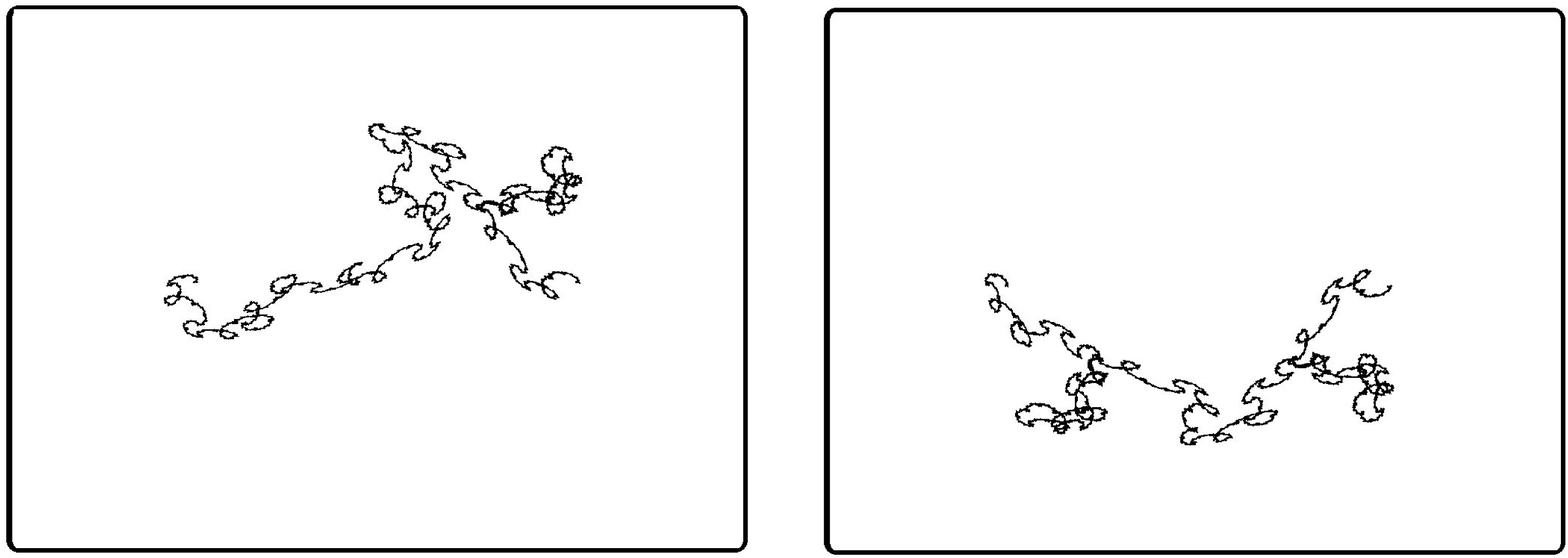}\\ 
$D(L,L) $  \hspace{1.9cm} $U(R,R) $ \hspace{2.6cm} $U(R,L) $   \hspace{1.9cm} $D(L,L)$
\end{center} 

\begin{center}
\includegraphics[scale=.09]{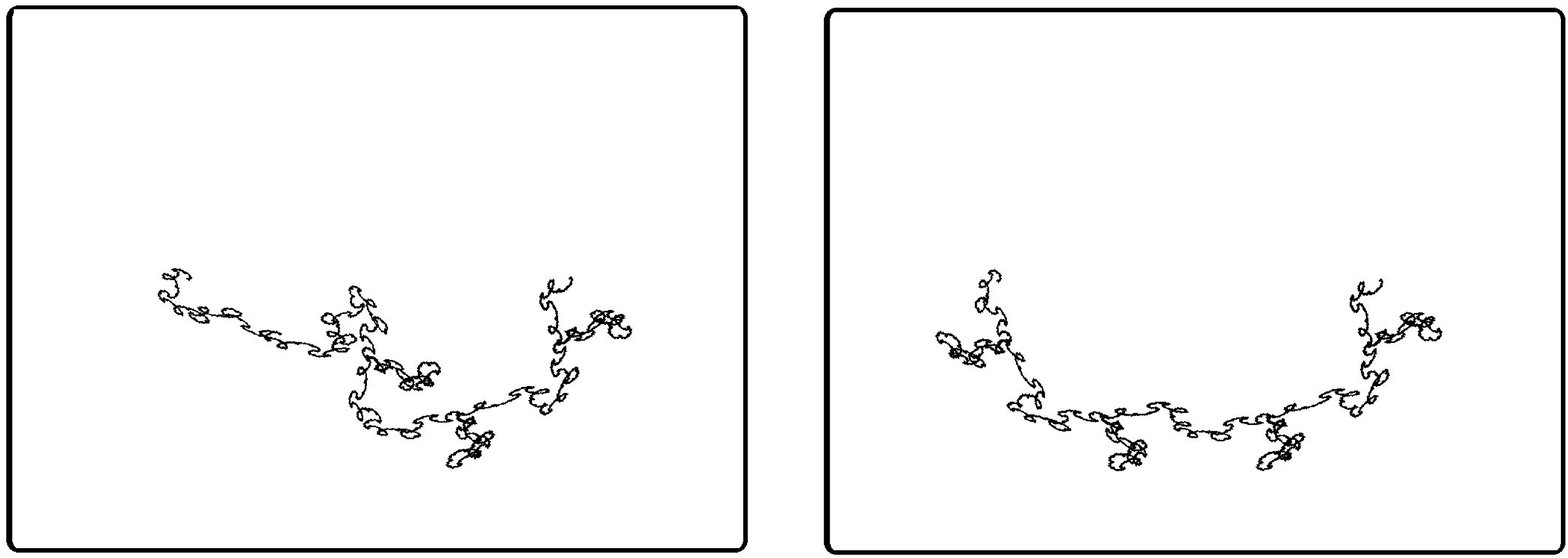}  \qquad  
\includegraphics[scale=.09]{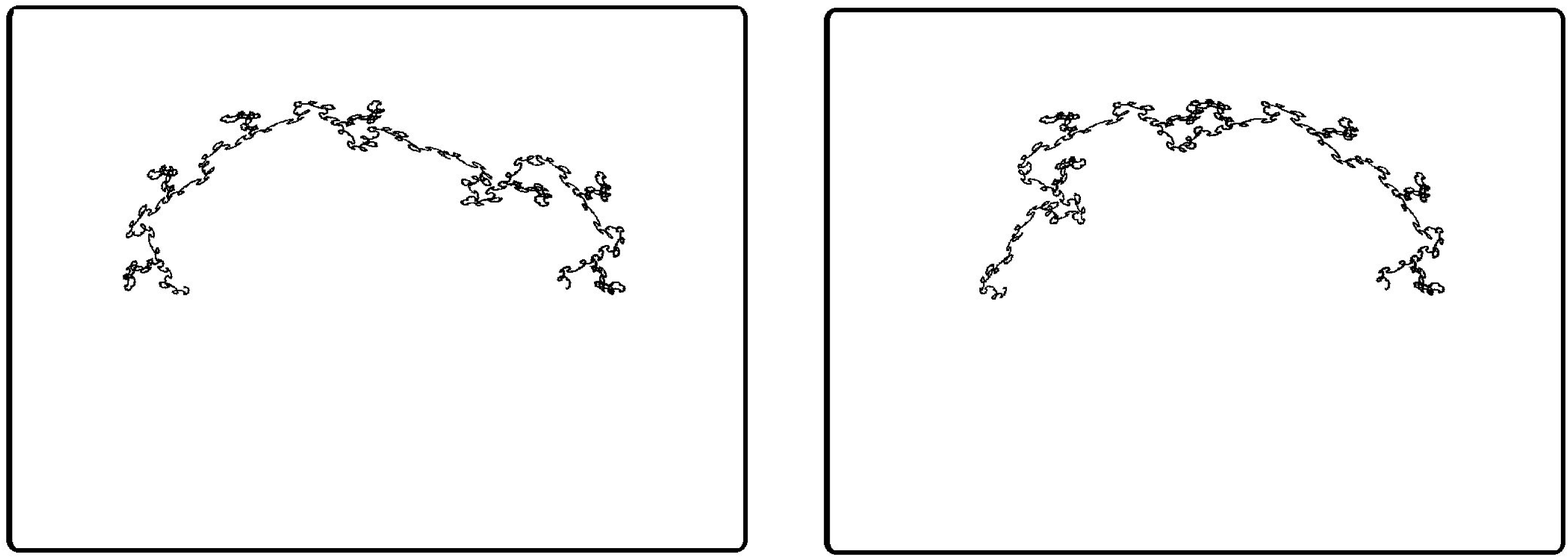}\\ 
$D(L,R) $  \hspace{1.9cm} $D(R,R) $ \hspace{2.6cm} $U(R,L) $   \hspace{1.9cm} $U(L,R)$
\end{center} 

\begin{center}
\includegraphics[scale=.09]{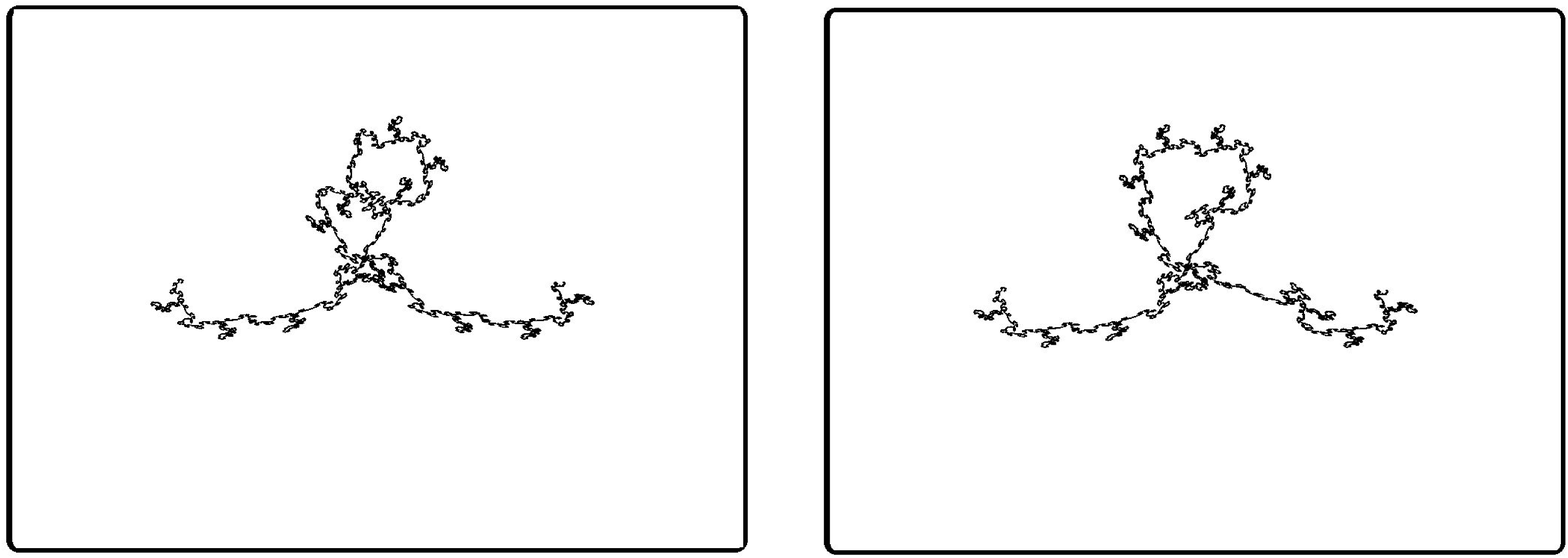}  \qquad
 \includegraphics[scale=.09]{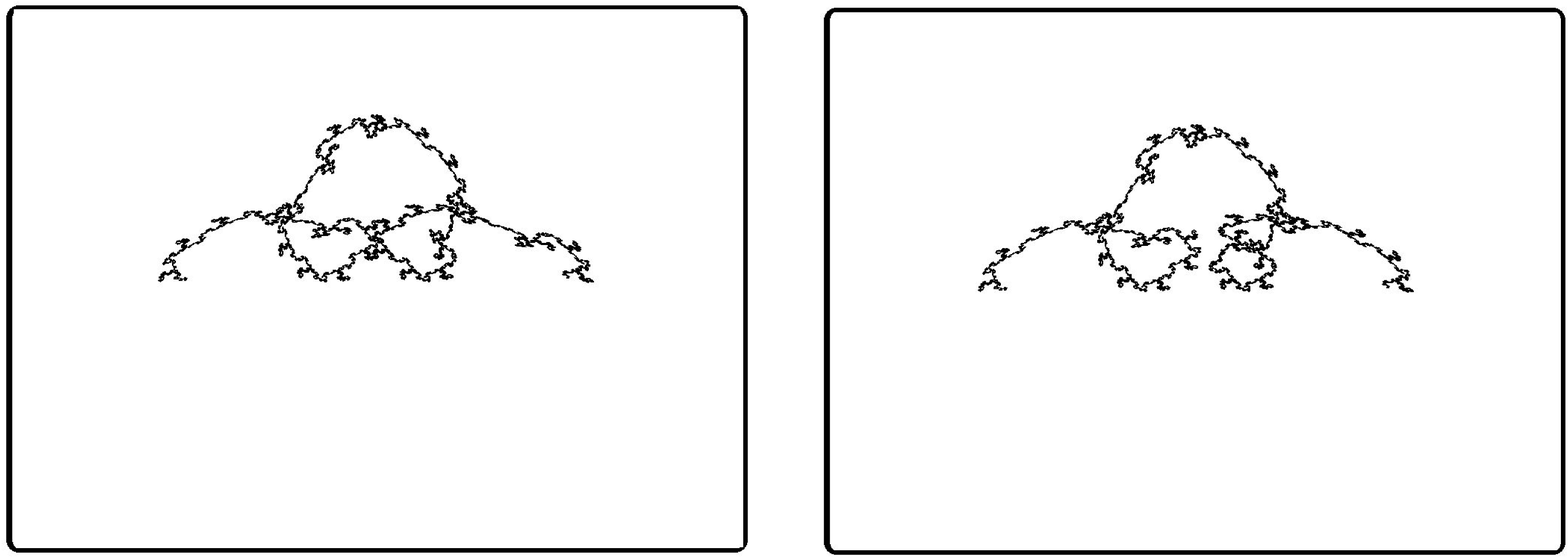}\\ 
$U(L,R) $  \hspace{1.9cm} $U(L,L) $ \hspace{2.6cm} $U(R,R) $   \hspace{1.9cm} $U(R,L)$
\end{center} 

\begin{quote}
\textsf{\textbf{Figure 6.} A sequence of pairs of (approximate) 2-variable fractals  (from left to right and then down the page).}
\end{quote} 

\medskip

In each case in this particular example, for each buffer at each step, either $ U $ or $ D $ was chosed with probability $\frac{1}{2} $.  For each buffer, the (in this case two) input buffers chosen from the previously generated  pair of buffers were also each chosen as either $ L $ or $ R $ with probability $\frac{1}{2} $, and the same buffer is allowed to be selected twice. 
After about 12 iterations, the images obtained are independent of the initial images up to screen resolution.  After this stage the images (or ``necklaces'') can be considered as examples of 2-variable fractals corresponding to the family $  (U,D) $ of IFSs with associated choice probabilities  $  (\frac{1}{2} ,\frac{1}{2} )  $.  

The pair of 2-variable fractals obtained at each step  depends on the previous choices of IFS and input buffers, and will vary from one experimental run to another.  However, over any sufficiently long experimental run, the empirically obtained distribution on pairs of 2-variable fractals will (up to any prescribed resolution) be the same with probability one.  This follows from ergodic theory and the fact that the construction process corresponds to the chaos game for an IFS (operating here on pairs of images rather than on single points as does a standard IFS).  As discussed in the previous section, we call this type of IFS a \emph{superIFS}.  The collection of 2-variable necklaces obtained over a long experimental run should be thought of as a single \emph{superfractal}, and the corresponding probability distribution on necklaces should be thought of as the corresponding \emph{superfractal measure}.  

In Figure 7 we have superimposed the members of a generated sequence of 2-variable fractal necklaces.  By virtue of the fact that, as discussed before,  the probability distribution given by such a  sequence  approximates the associated superfractal measure, the image can be regarded as a projection of the superfractal onto 2-dimensional space.  The attractors of the individual IFSs $U$ and $D$ are shown in green and red respectively.  The projected support of the superfractal is shown on a black background but, inside the support, increasing density of the superfractal measure is indicated by increasing intensity of white.

\begin{center}
\includegraphics[scale=.3]{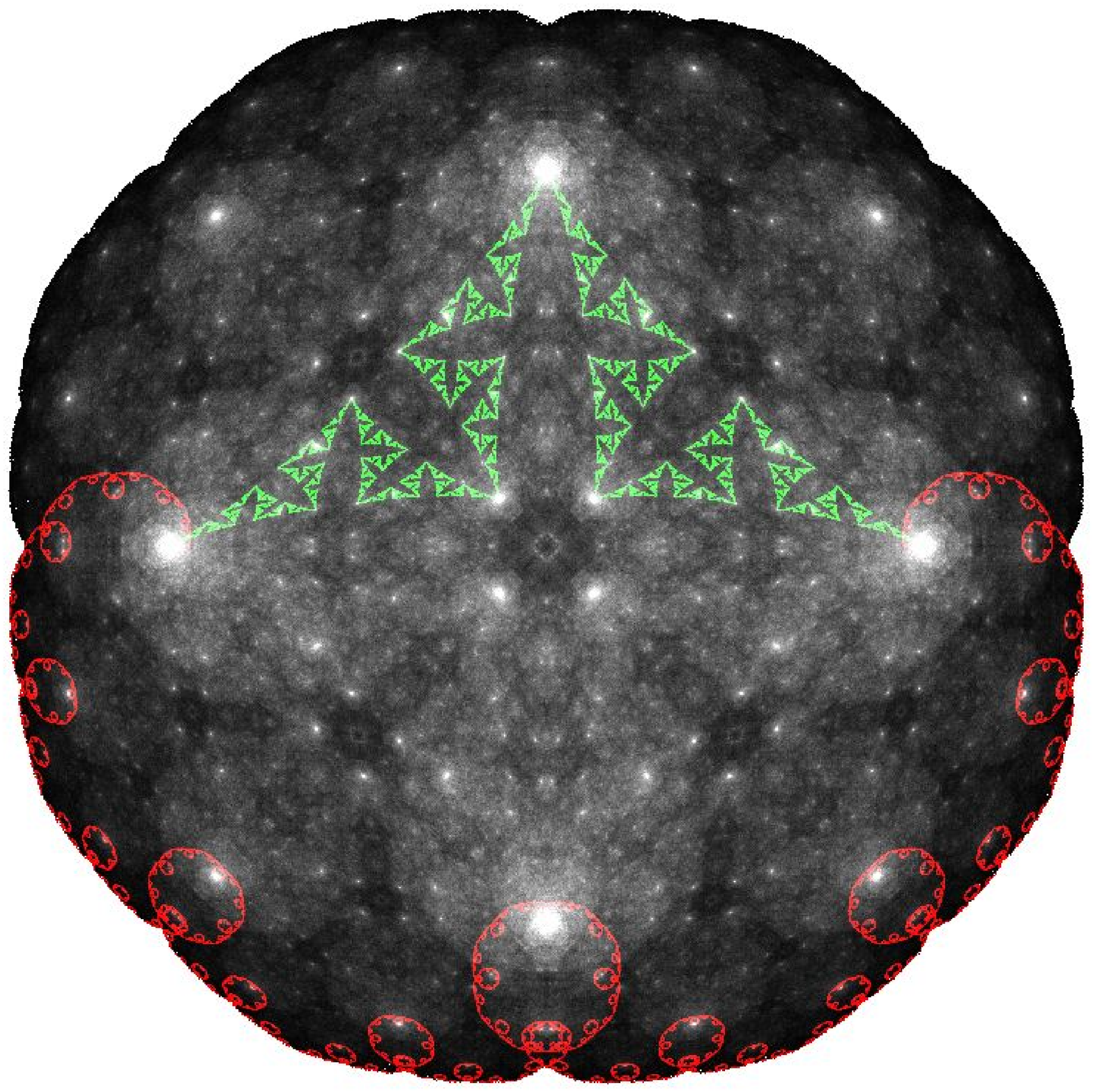}\\
\textsf{\textbf{Figure  7.}  Superfractal projected onto 2 dimensional
 space.}
 \end{center}

\medskip 
In Figure 8 are shown a fern and lettuce generated by two IFSs, each IFS consisting of 4 functions.  In Figure 9  is a sequence of hybrid offspring, extending the examples in 
Figure  4.  
The colouring was obtained by working with two IFSs in 5 dimensional space, with the three additional dimensions corresponding to RGB colouring. The two IFSs used project onto two IFSs operating in two dimensional space and which give the (standard black and white) fern and  lettuce attractors respectively.  The 2-variable offspring were coloured by extending the superfractal construction to 5 dimensional space in a natural manner.

\begin{center}
\includegraphics[scale=.5]{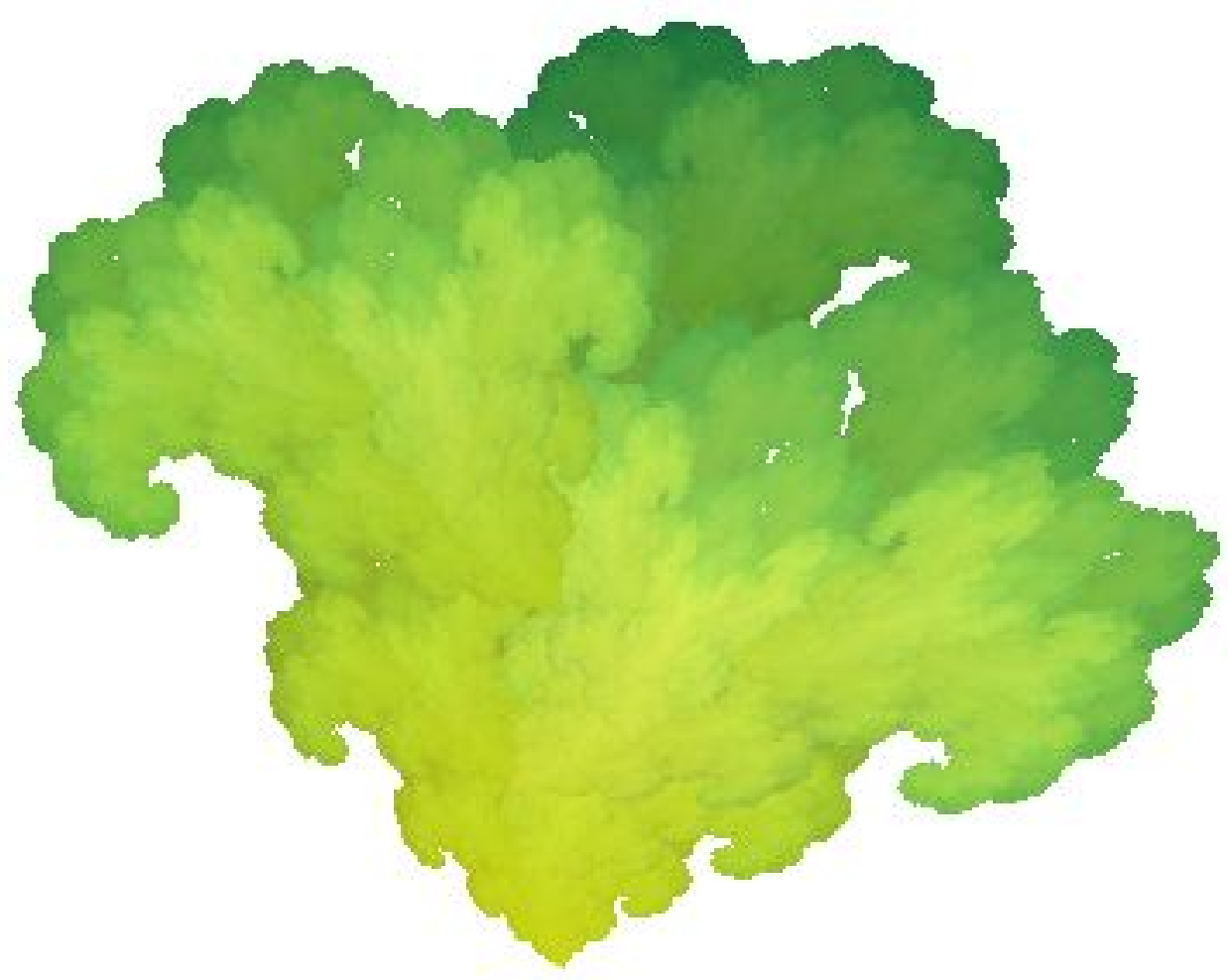} \quad  
\includegraphics[scale=.5]{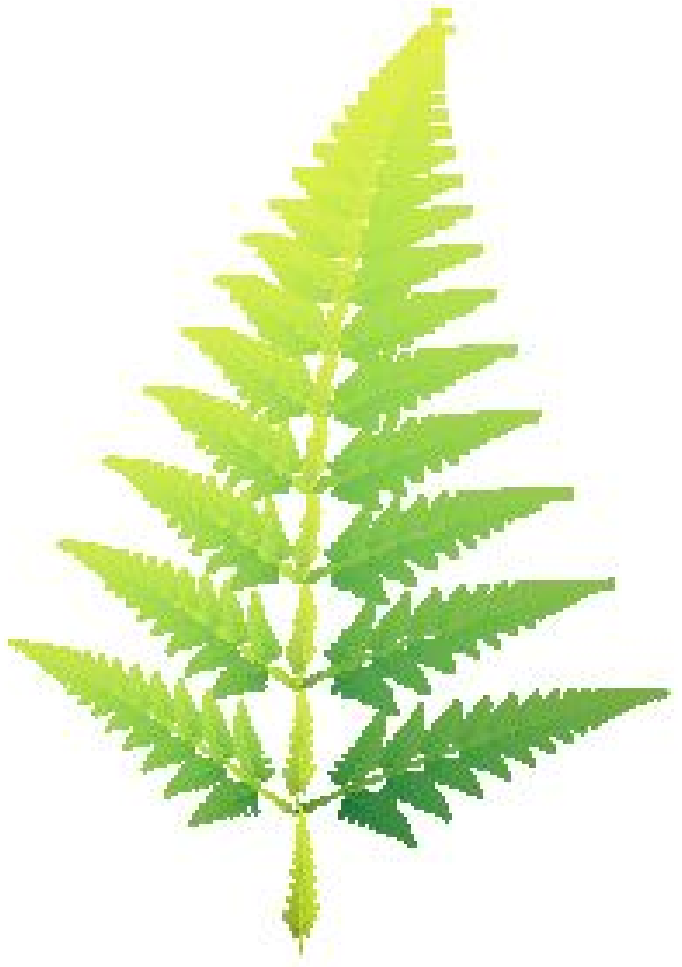} 

\textsf{\textbf{Figure  8.}  Lettuce and fern attractors.}
\end{center} 

\medskip

\begin{center}
\includegraphics[scale=.2]{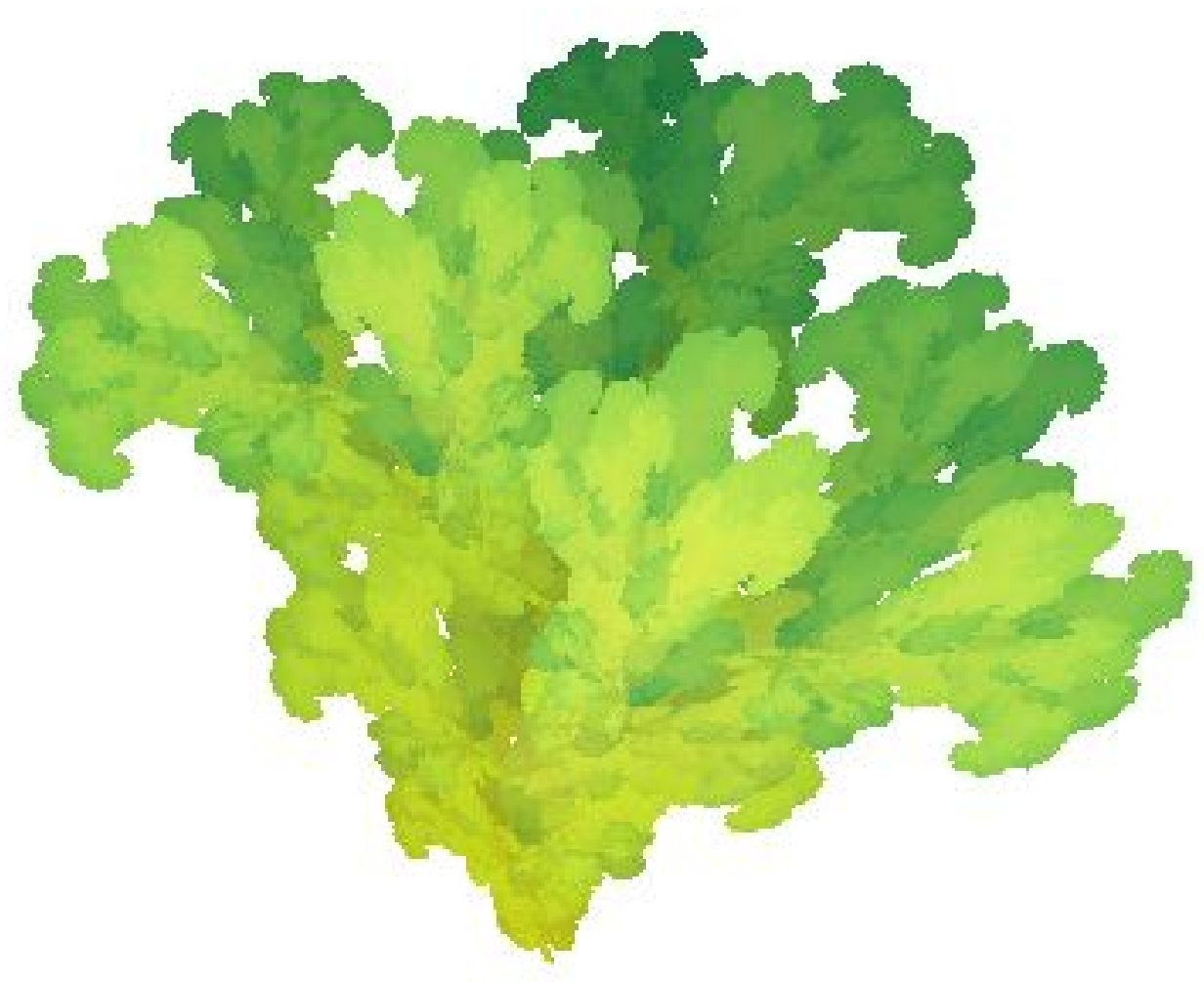} \quad  
\includegraphics[scale=.2]{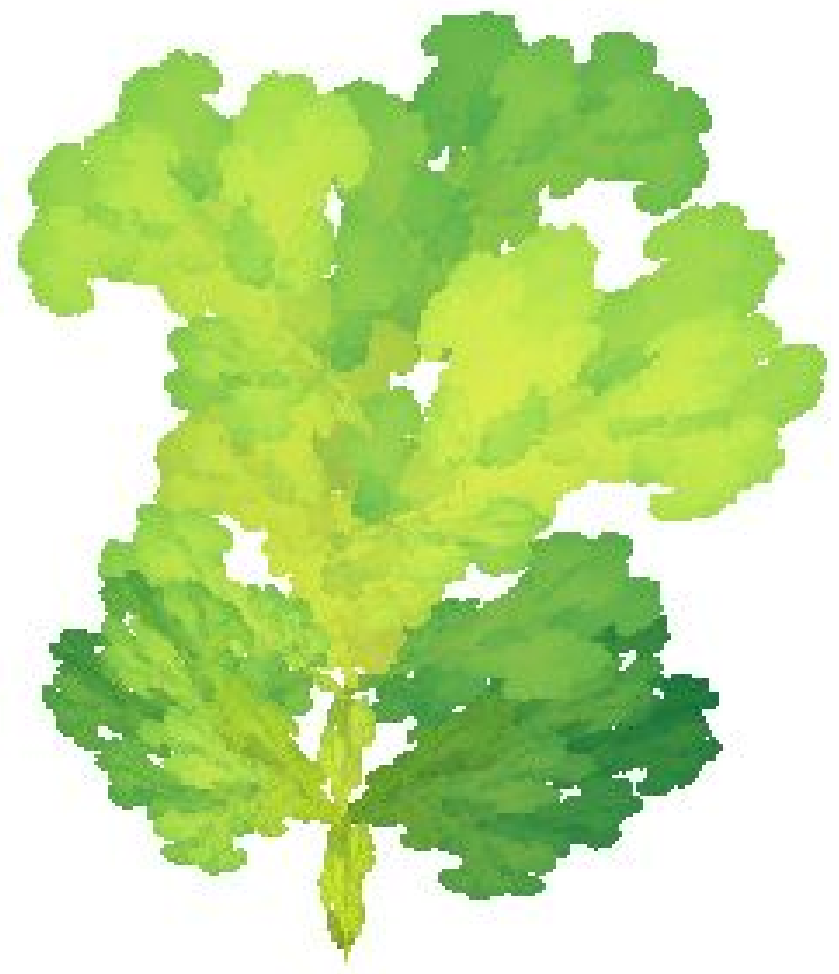} \quad 
\includegraphics[scale=.2]{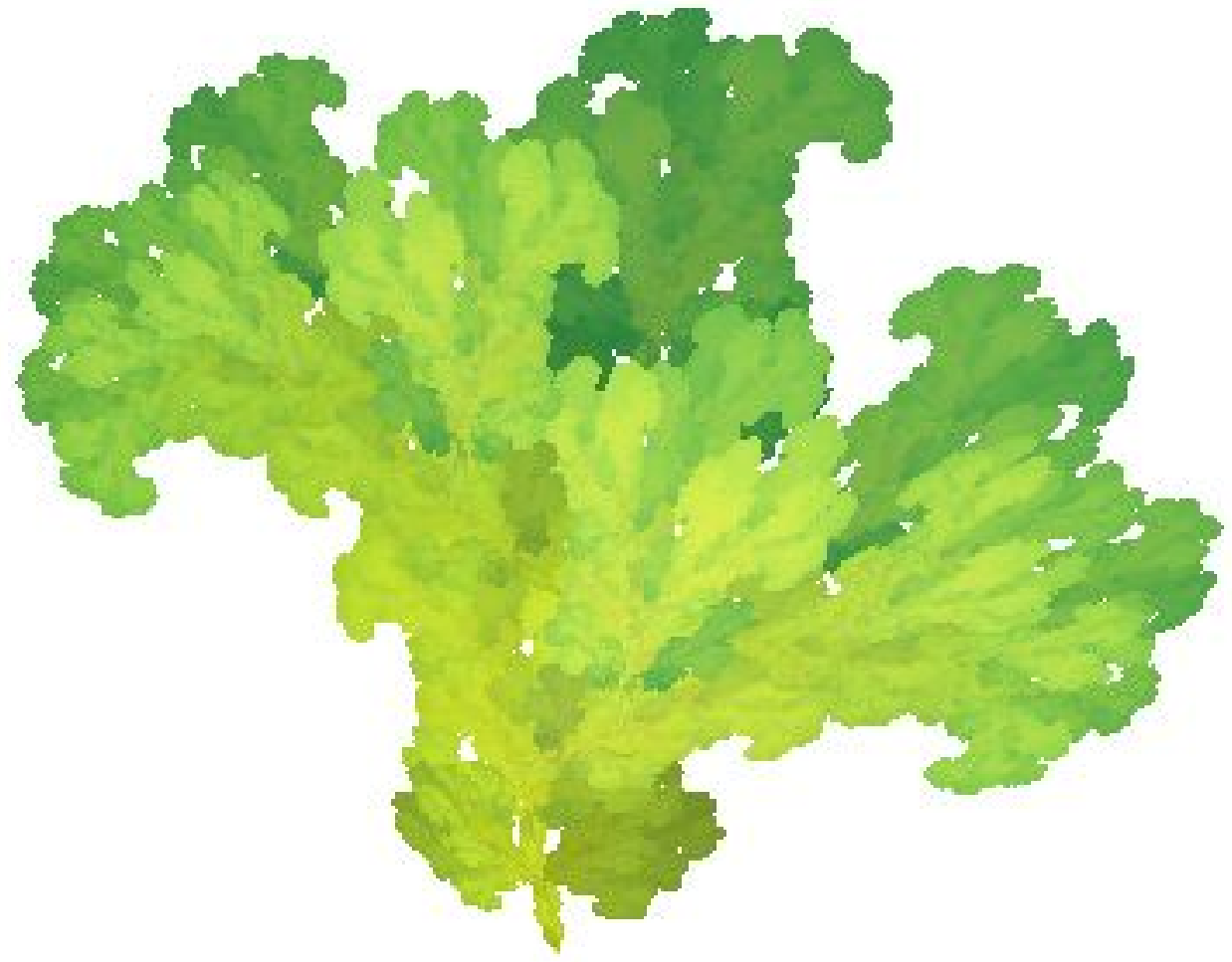}  \quad  
\includegraphics[scale=.2]{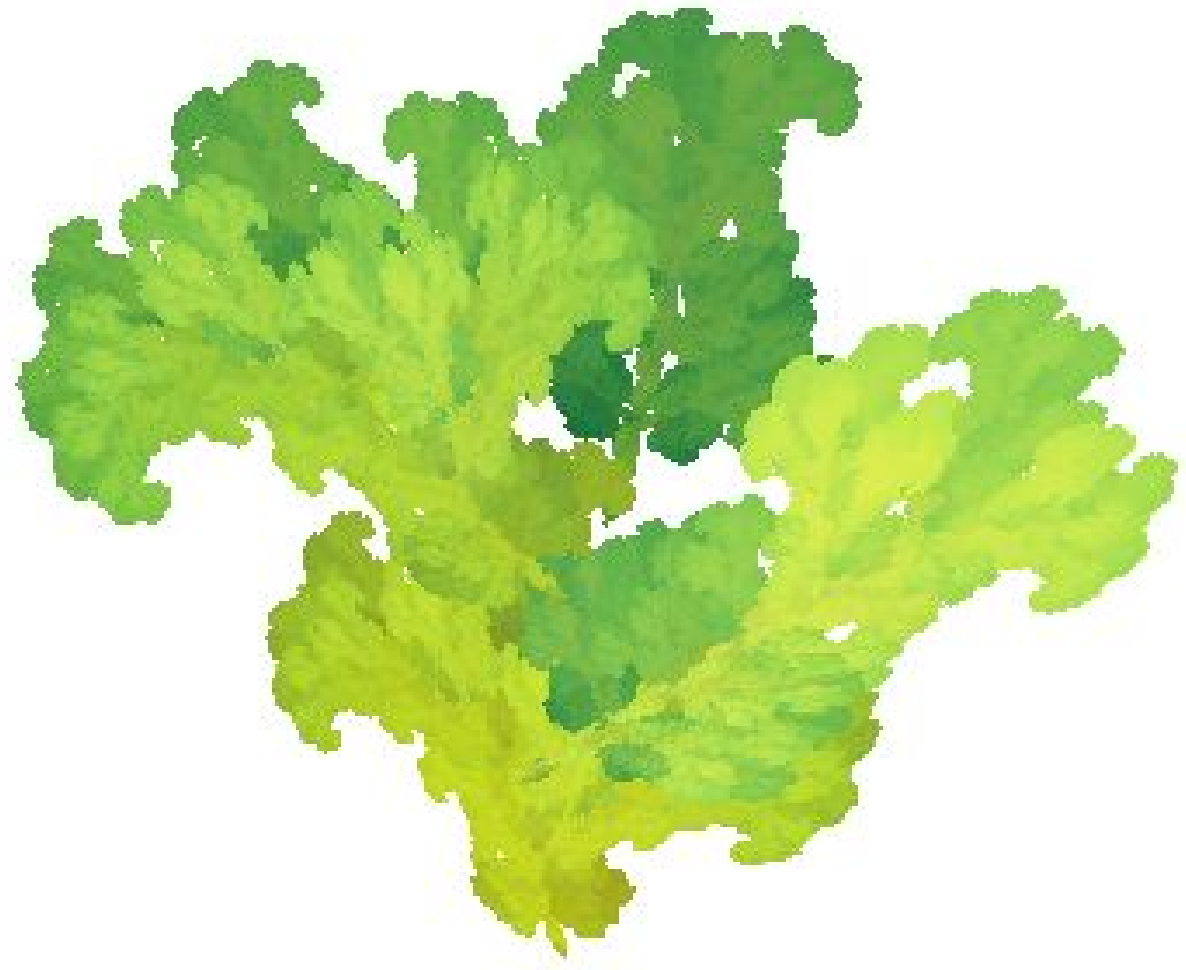} \quad 
\includegraphics[scale=.2]{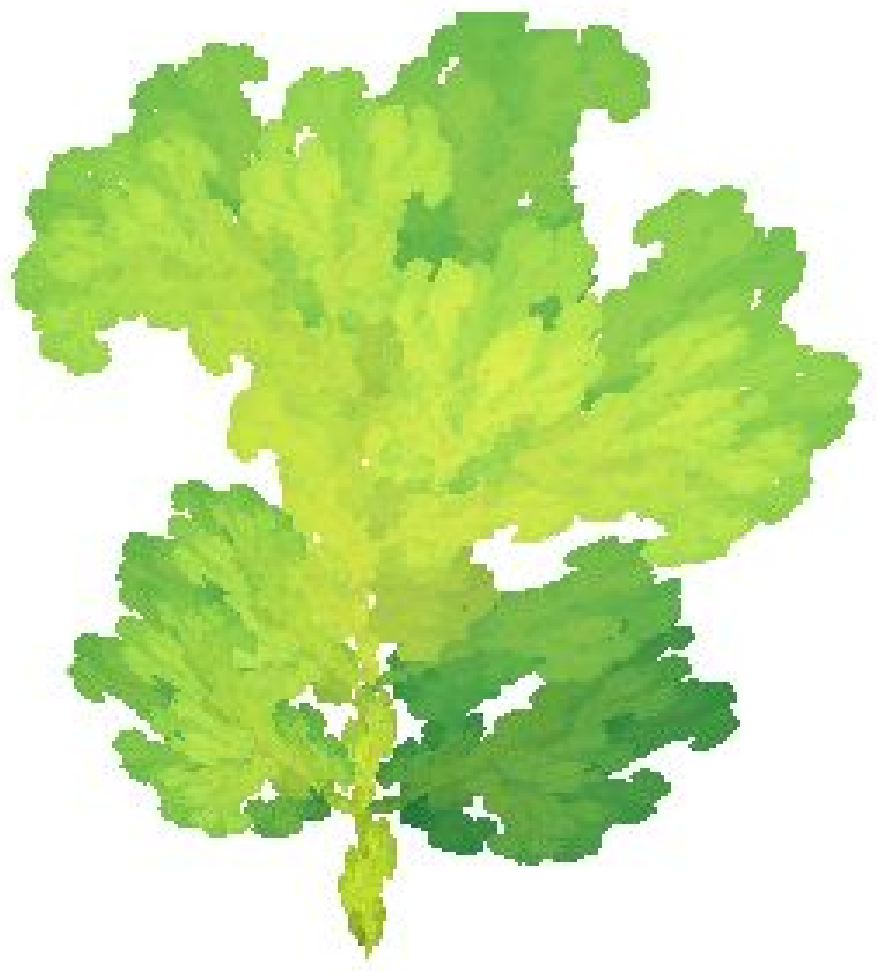} \quad 
\includegraphics[scale=.2]{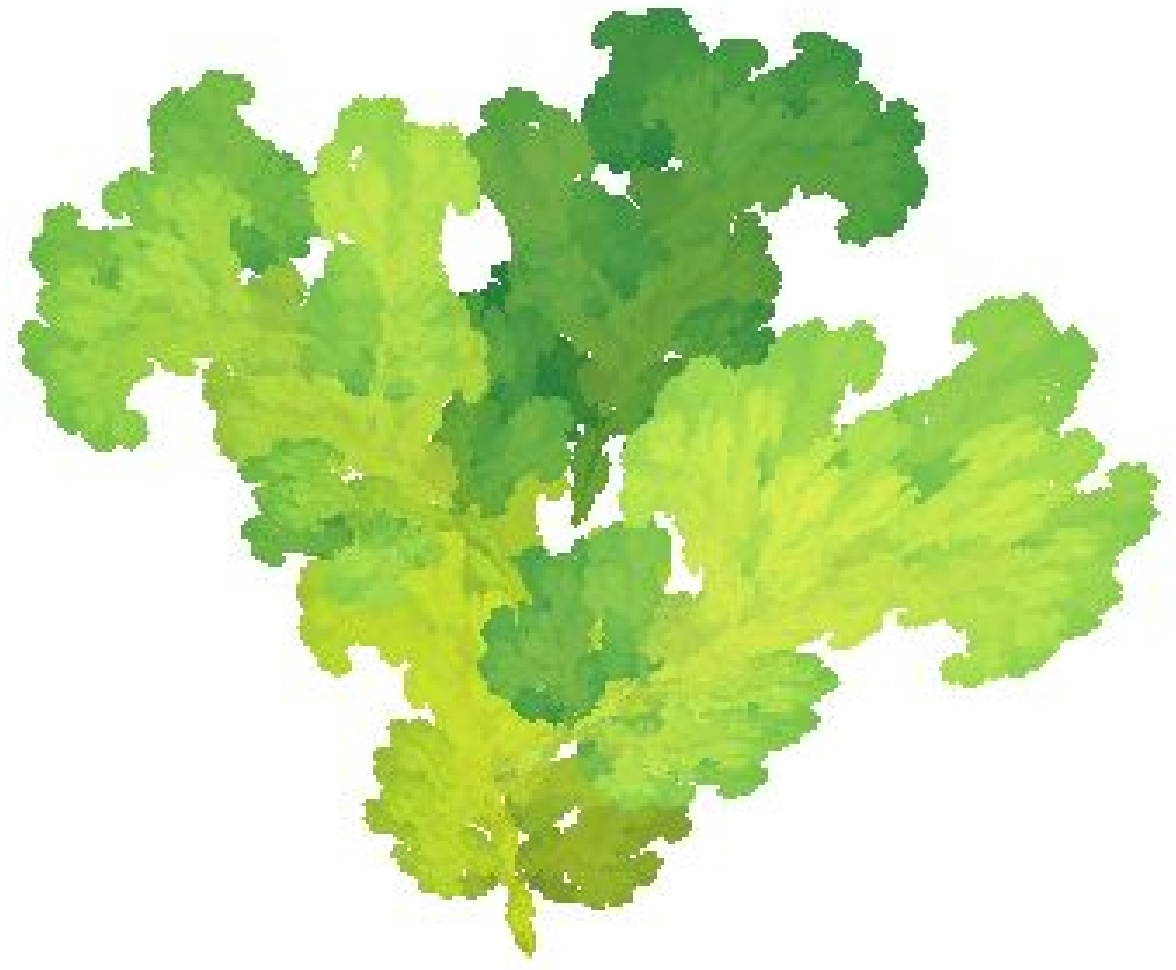}  \quad  
\includegraphics[scale=.2]{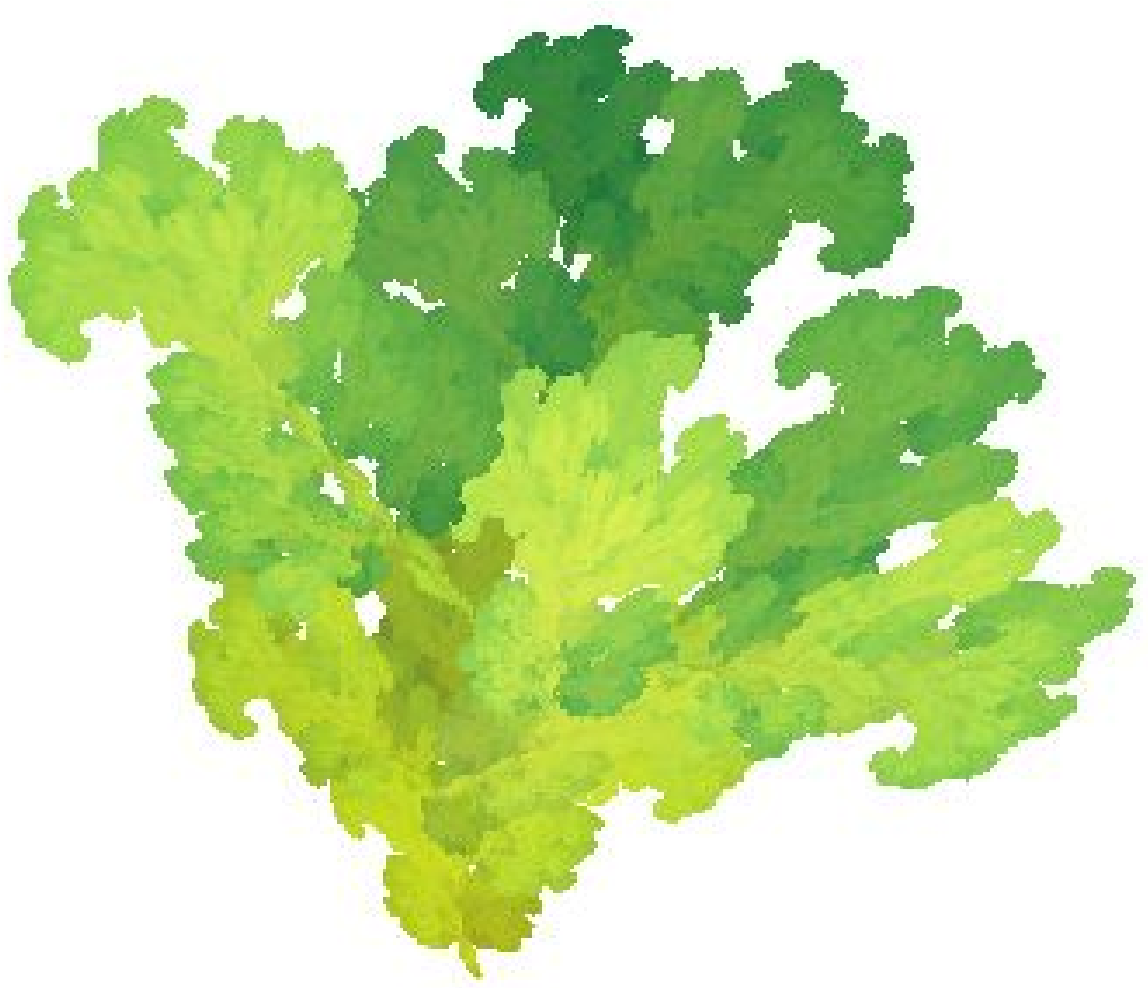} \quad 
\includegraphics[scale=.2]{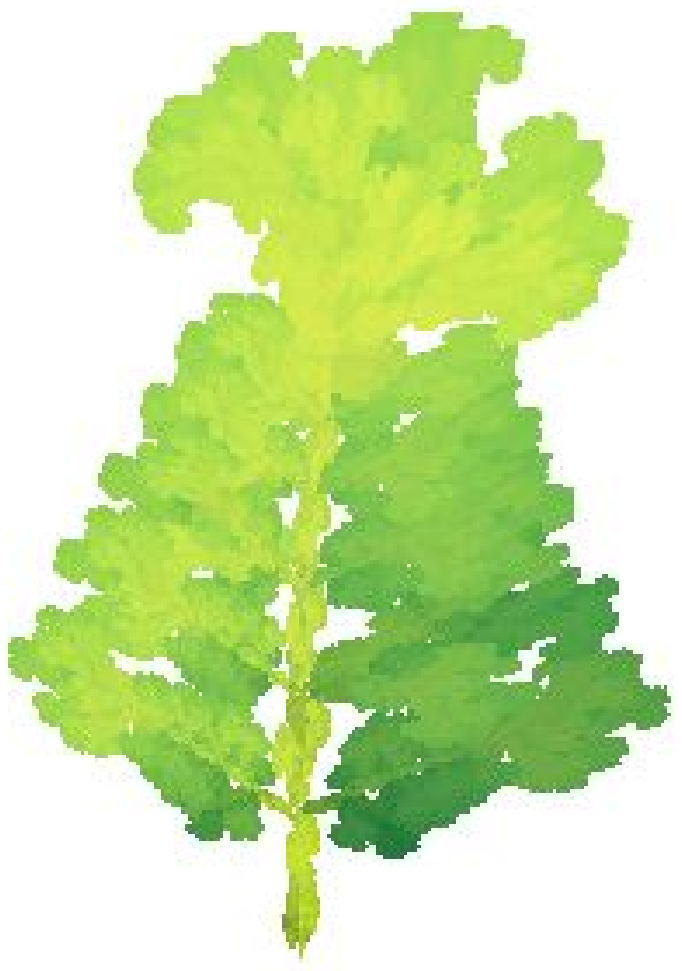} \quad  
\includegraphics[scale=.2]{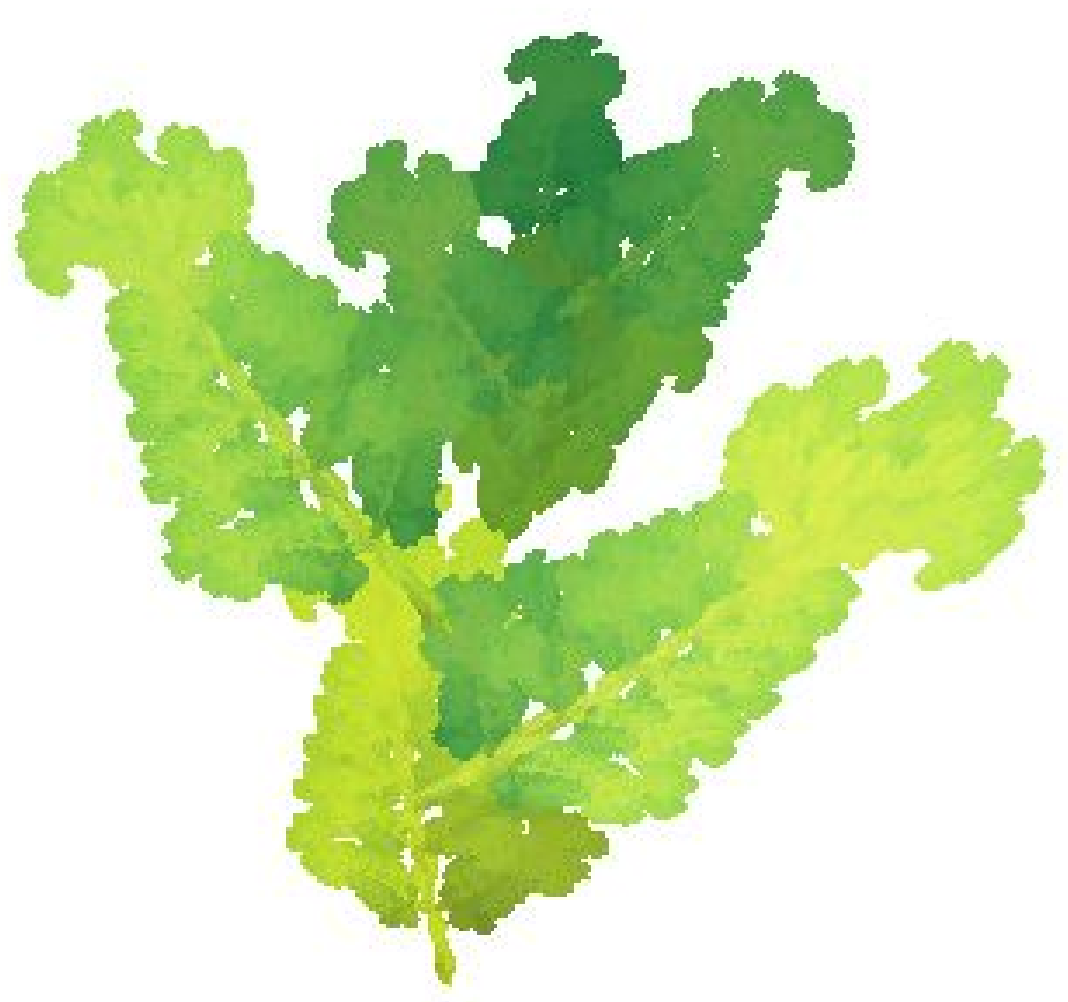} \quad 
\includegraphics[scale=.2]{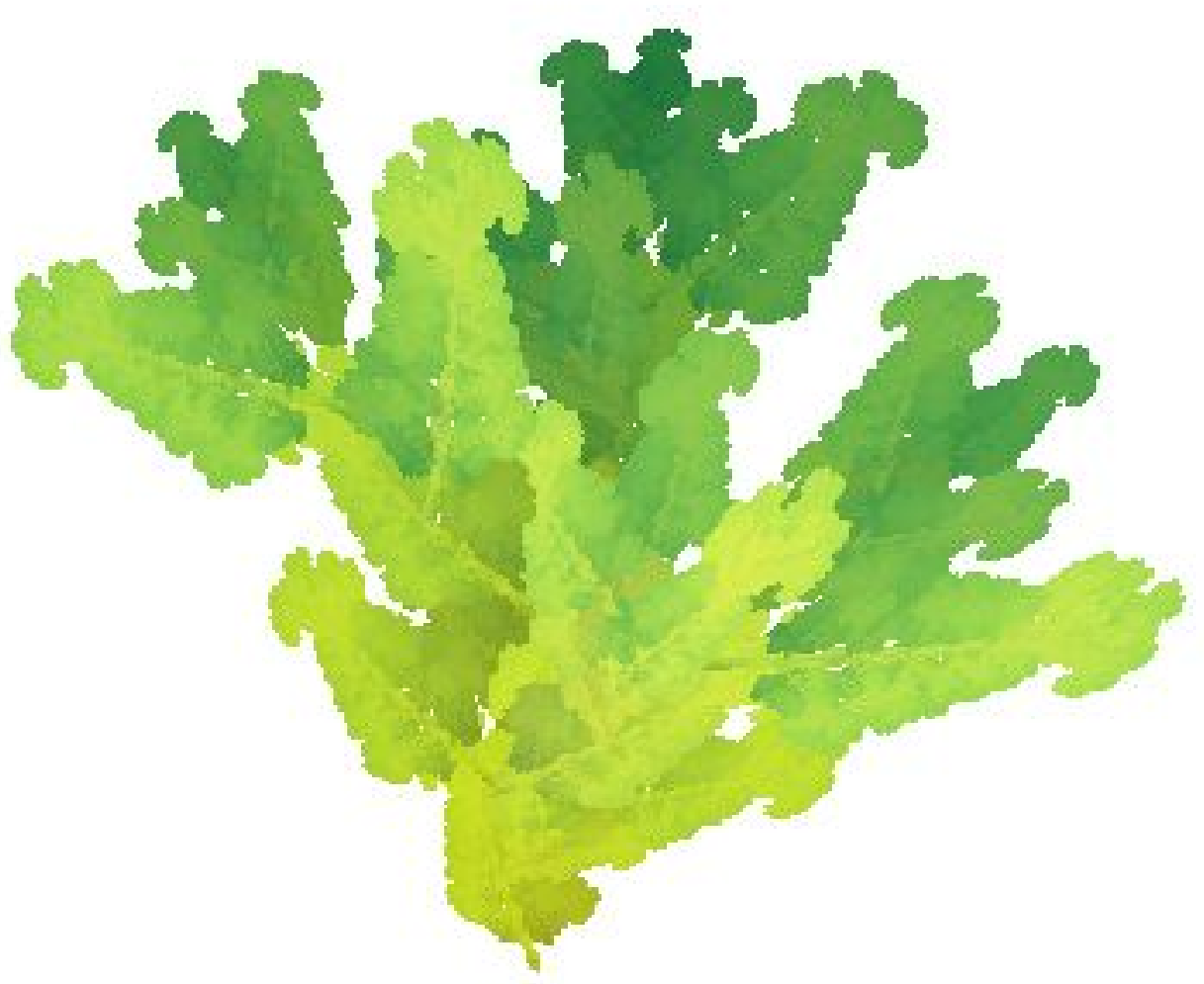}  \quad  
\includegraphics[scale=.2]{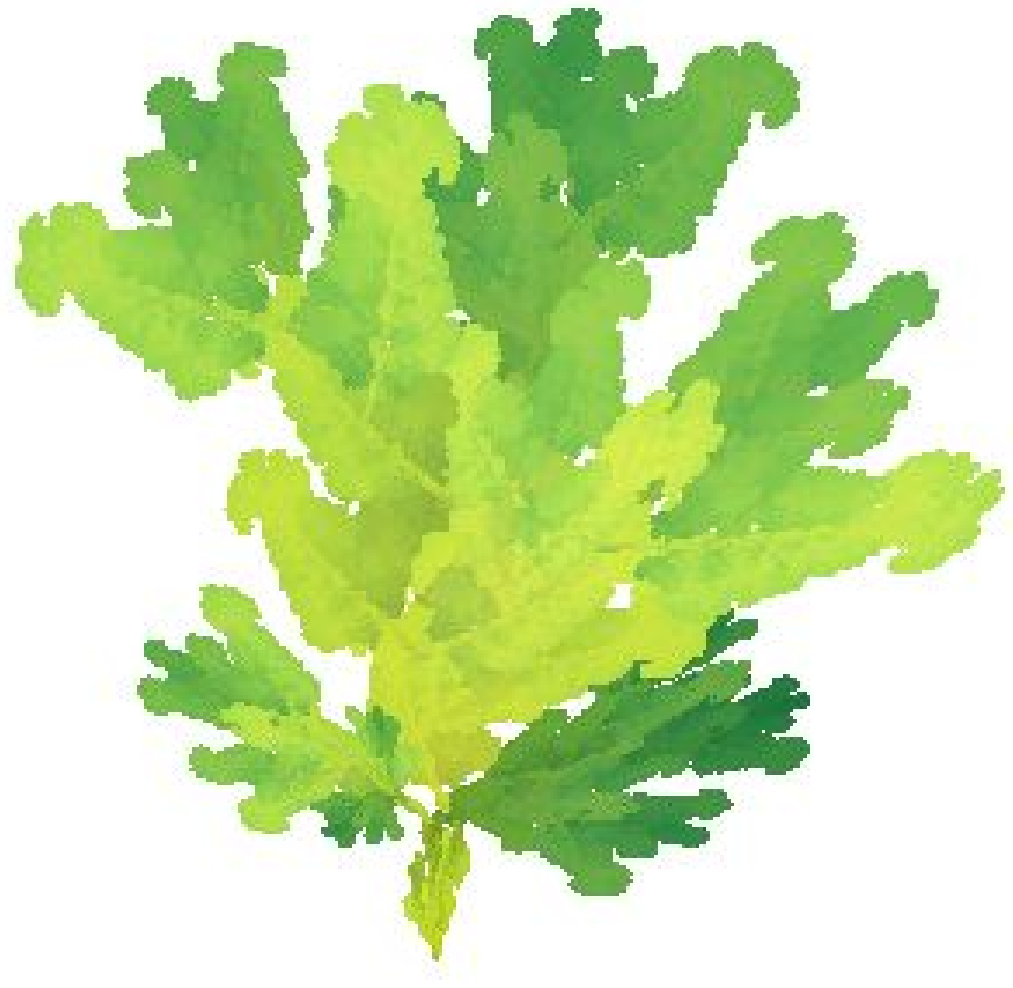} \quad 
\includegraphics[scale=.2]{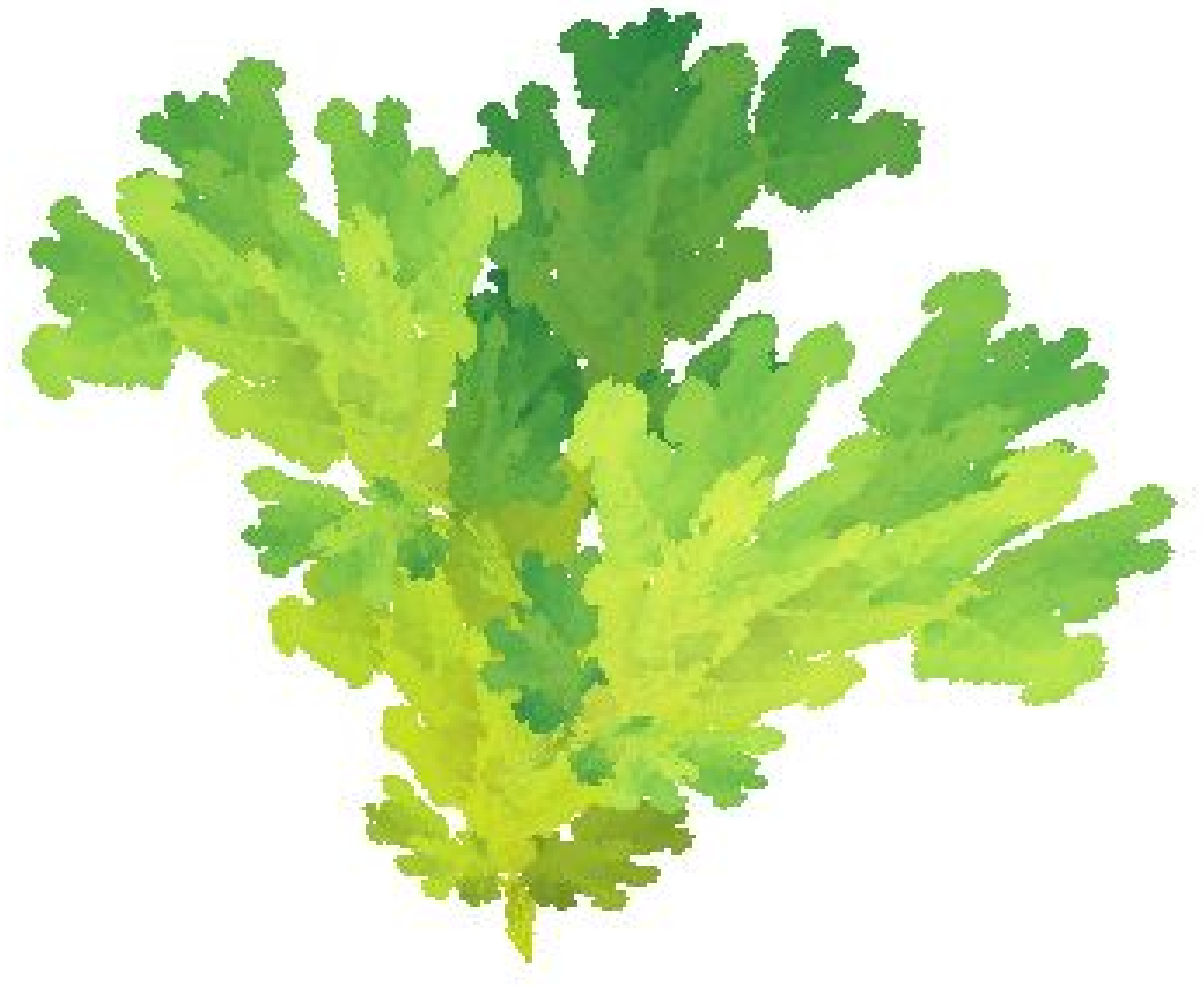} \quad  
\includegraphics[scale=.2]{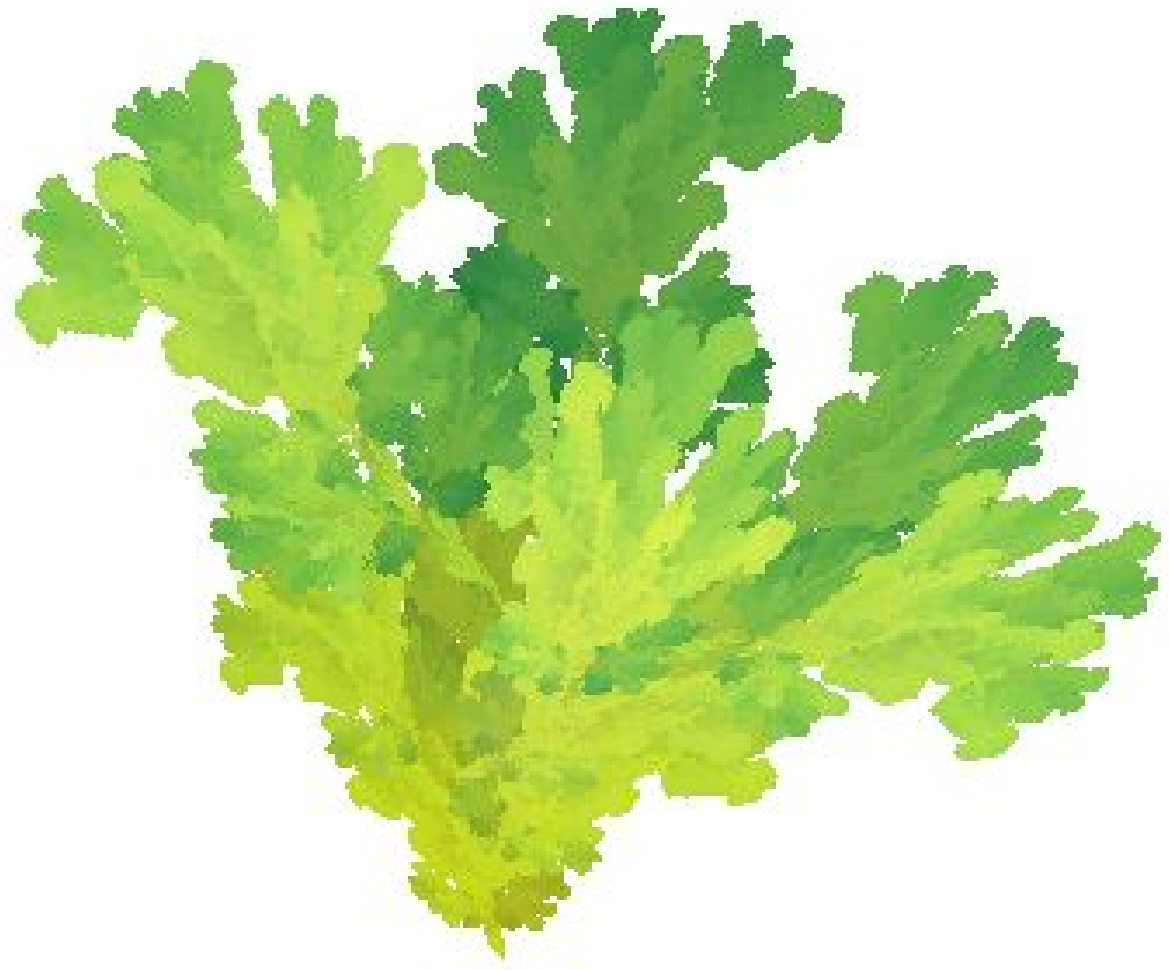} \quad 
\includegraphics[scale=.2]{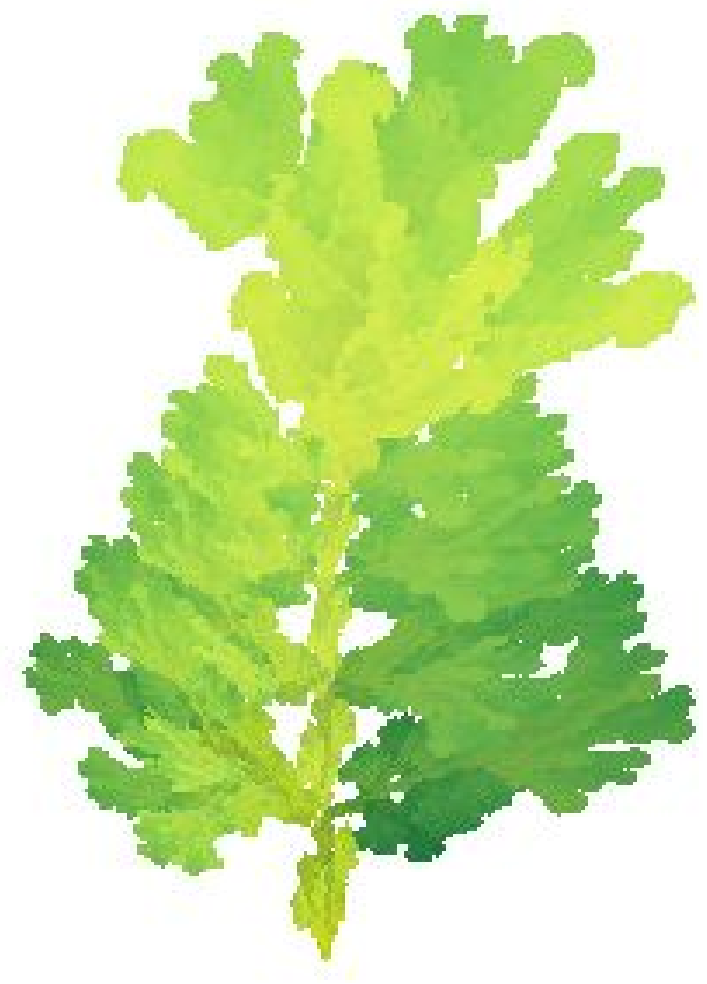}  \quad  
\includegraphics[scale=.2]{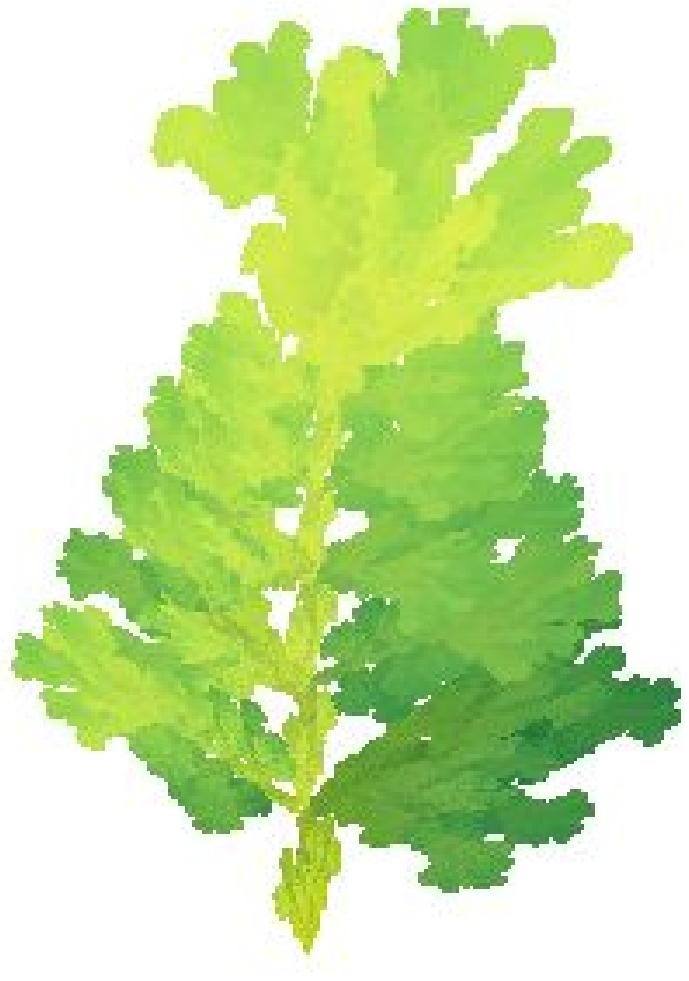} \quad 
\includegraphics[scale=.2]{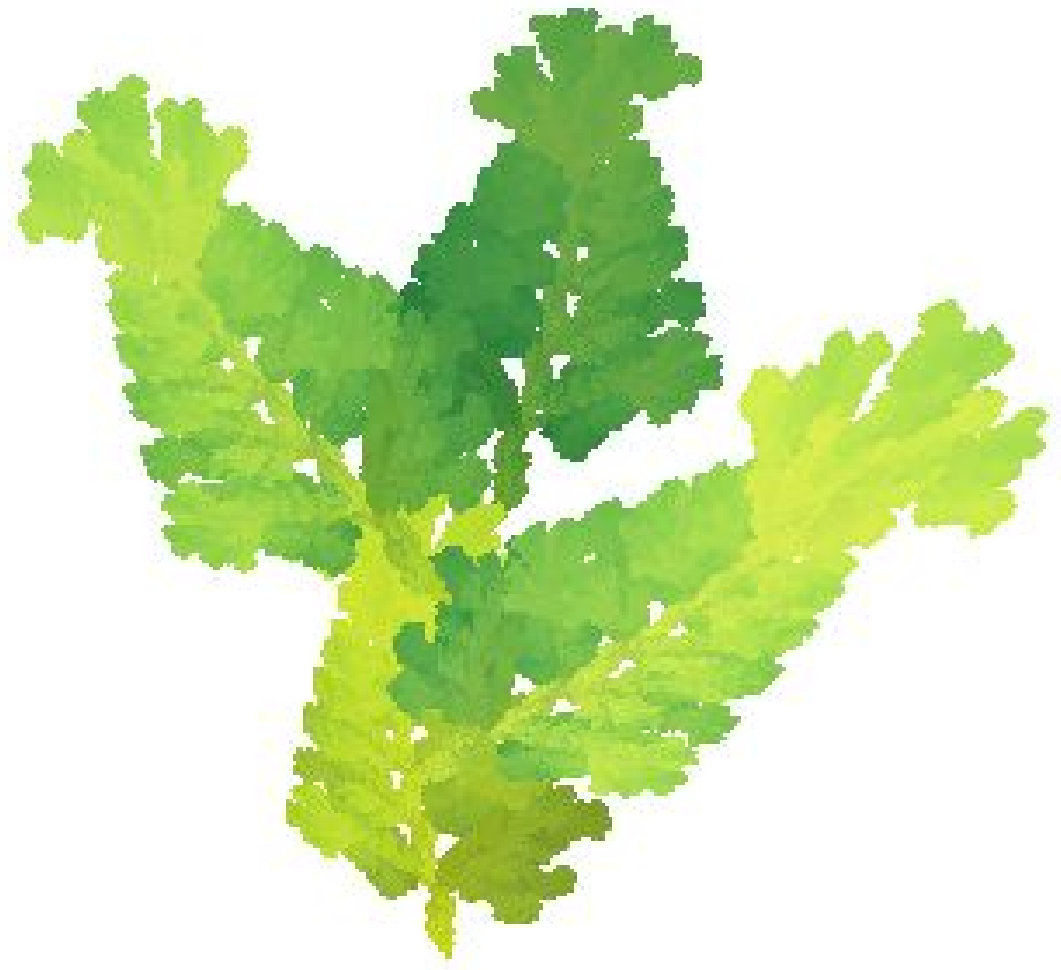} \quad  
\includegraphics[scale=.2]{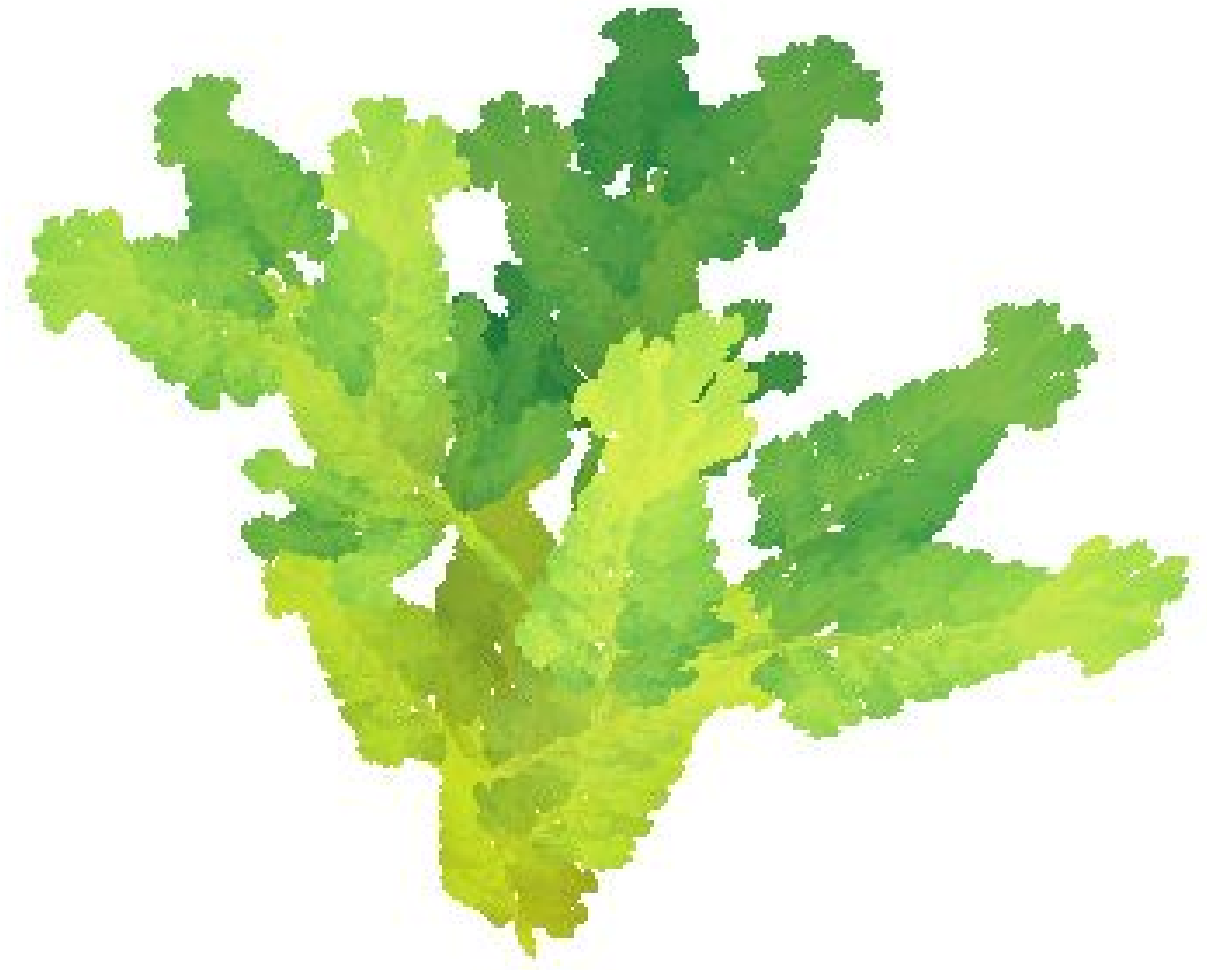} \quad 
\includegraphics[scale=.2]{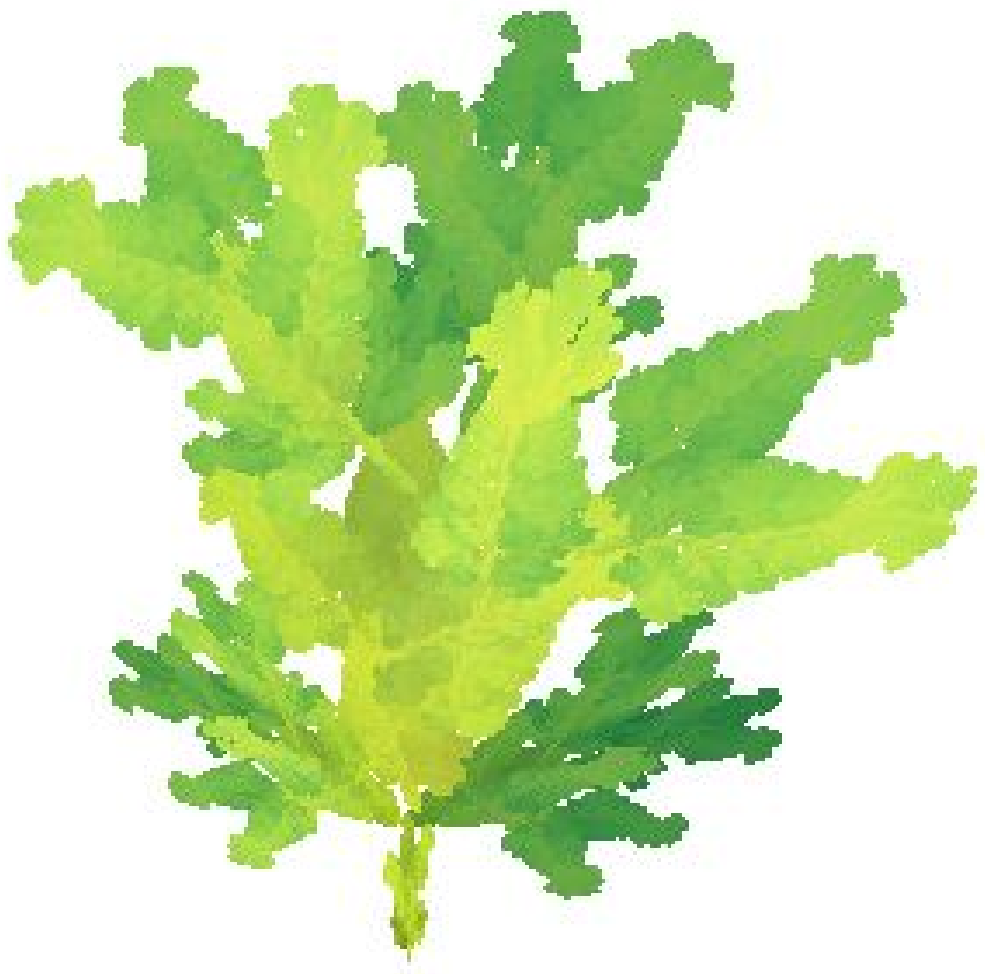}  \quad  
\includegraphics[scale=.2]{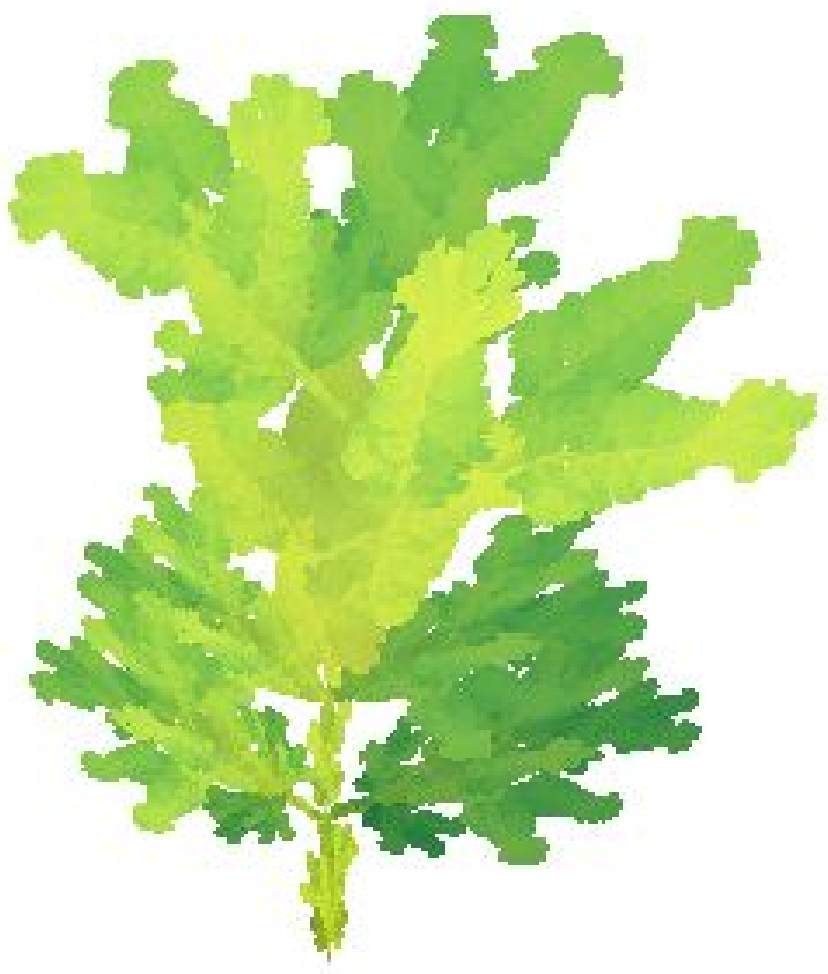} \quad 
\includegraphics[scale=.2]{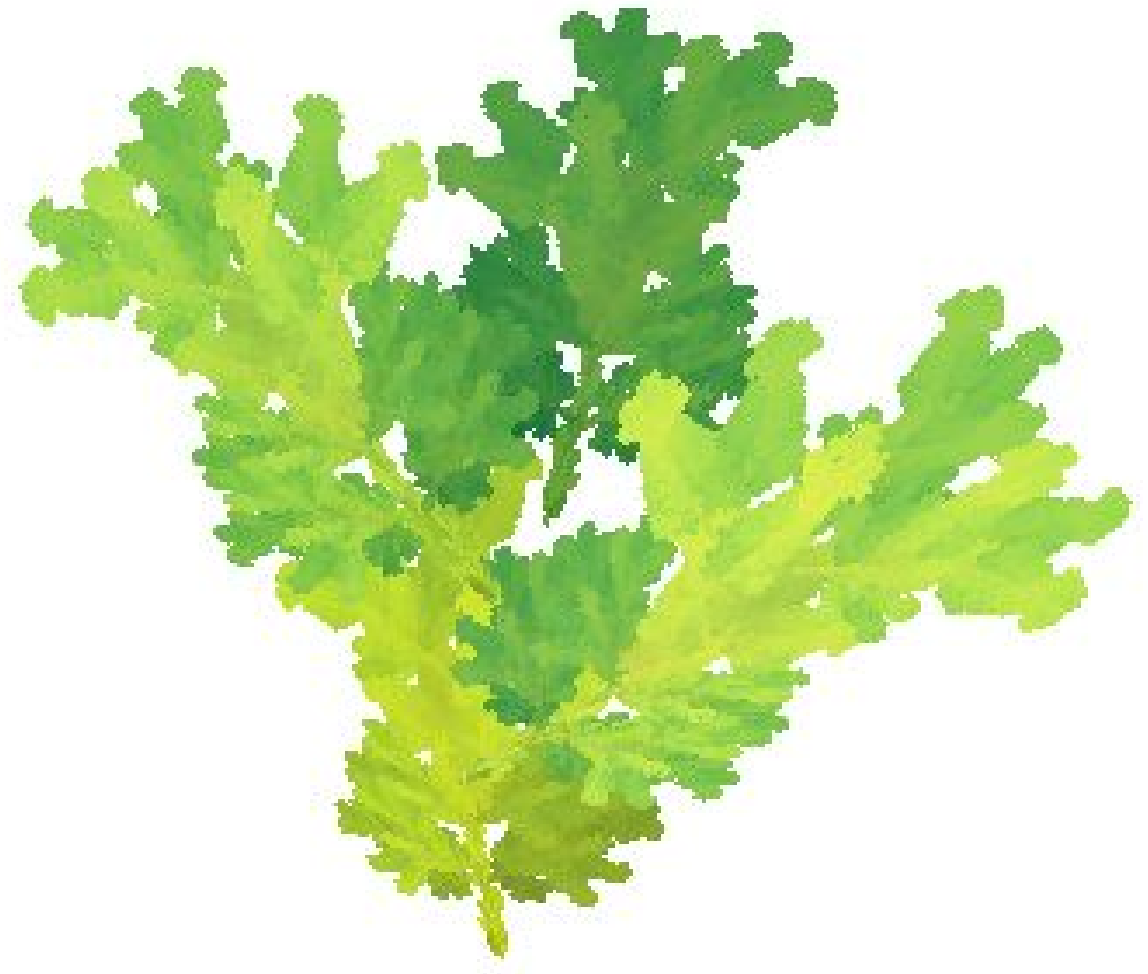} \quad\\ 
\textsf{\textbf{Figure 9.} A sequence of fern-lettuce hybrid offspring.}
\end{center}


\section{Computation of Dimension for $ V $-variable Fractals}
An important theoretical and empirical classification of fractals is via their dimension.  We first show how to compute the dimensions of  $V$-variable Sierpinski triangles.  As we discuss in the next section, the method generalises.

Associated with the transition from the  $k$-level set of $ V $ buffers to the   $ k+1 $ level is a  matrix  $ M^k(\alpha) $  defined for each  $ \alpha  $  as follows.  Entries of  $ M^k(\alpha) $ are initialised to zero.  For the set in the $ v $th output buffer at level  $ k+1 $  one considers each input buffer $ w $   used in its construction and adds  $ r^\alpha  $  to the $ w $th entry in the  $ v $th row where $r=\frac{1}{2}  $ or $  \frac{1}{3}   $   is the corresponding contraction ratio.   The construction of $ M^k(\alpha) $ can be seen in passing from level 1 to level 2 and from level 2 to level 3 in Figure 3, giving respectively:
$$
 M^1(\alpha)=
\begin{bmatrix}
      &  \frac{1}{3^\alpha}   & \frac{1}{3^\alpha}     
      & \frac{1}{3^\alpha}   & 0          &   0     \\
      &  0   &   \frac{2}{2^\alpha} 
      & 0  & 0          &   \frac{1}{2^\alpha}      \\
      &  0   & \frac{1}{3^\alpha}     
      & \frac{1}{3^\alpha}    & 0          &   \frac{1}{3^\alpha}       \\
      &  \frac{1}{2^\alpha}   &   0   
      & 0 &  \frac{2}{2^\alpha}       &   0     \\
      & 0   &  0
      & \frac{1}{2^\alpha}    &\frac{1}{2^\alpha}         &  \frac{1}{2^\alpha}       \\ 
\end{bmatrix},
\qquad 
M^2(\alpha) = 
\begin{bmatrix}
      &  \frac{1}{2^\alpha}   & \frac{1}{2^\alpha}   
      & \frac{1}{2^\alpha}   & 0          &   0     \\
      &  \frac{1}{3^\alpha}    &   0
      & \frac{1}{3^\alpha}   & \frac{1}{3^\alpha}           &   0     \\
      & 0  &    0
      & \frac{1}{3^\alpha}    & \frac{1}{3^\alpha}          &    \frac{1}{3^\alpha}       \\
      &  \frac{1}{3^\alpha}   &   0   
      & 0 & 0   &  \frac{2}{3^\alpha}      \\
      & 0   &  0
      & \frac{2}{2^\alpha}    &0        &  \frac{1}{2^\alpha}       \\ 
\end{bmatrix}.
$$

The ``pressure'' function  
\begin{equation} \label{pf} 
\gamma_V(\alpha)=\lim_{k\to \infty} \frac{1}{k}
     \log\left( \frac{1}{V}\left\| M^1(\alpha)\cdot\ldots\cdot M^k(\alpha)\right\|\right)
\end{equation} 
 exists and is independent of the experimental run with probability one by a result of Furstenberg and Kesten (1960), see also Cohen (1988) for  the version required here. (By $ \|A\| $ we mean the sum of the absolute values of all entries in the matrix $ A $.)   The factor $ 1/V $ is not necessary in the limit, but is the correct theoretical and numerical normalisation, as we see in the next section.   (See Feng and Lau (2002) for another use of  Furstenberg and Kesten type results for computing dimensions of random fractals.)  
In case $ V=1$,
\begin{equation} \label{pf1}  
 \gamma_1(\alpha) = \left(1-\frac{\alpha}{2}\right)\log 3 - \frac{\alpha}{2}\log 2 
\end{equation} 
  from the strong law of large numbers.   
 
 It can be shown, see Barnsley, Hutchinson and Stenflo (2003b), that for each $V$ $\gamma_V( \alpha ) $ is monotone decreasing.  In this example   the derivative lies between $ -\log 2 $ and $- \log 3 $, corresponding to the contraction ratios $ \frac{1}{2} $  and $ \frac{1}{3}$   respectively, see Figure 10.
Moreover,  there is a unique $ d = d(V) $   such that $\gamma_V( d)= 0$.  This   is the dimension of the corresponding  $  V $-variable random fractals with probability one.  The establishment and generalisation of this method uses the theory of products of random matrices and ideas from statistical mechanics, as we discuss in the next section.  

\begin{center}
\noindent\includegraphics[scale=.6, angle=-90]{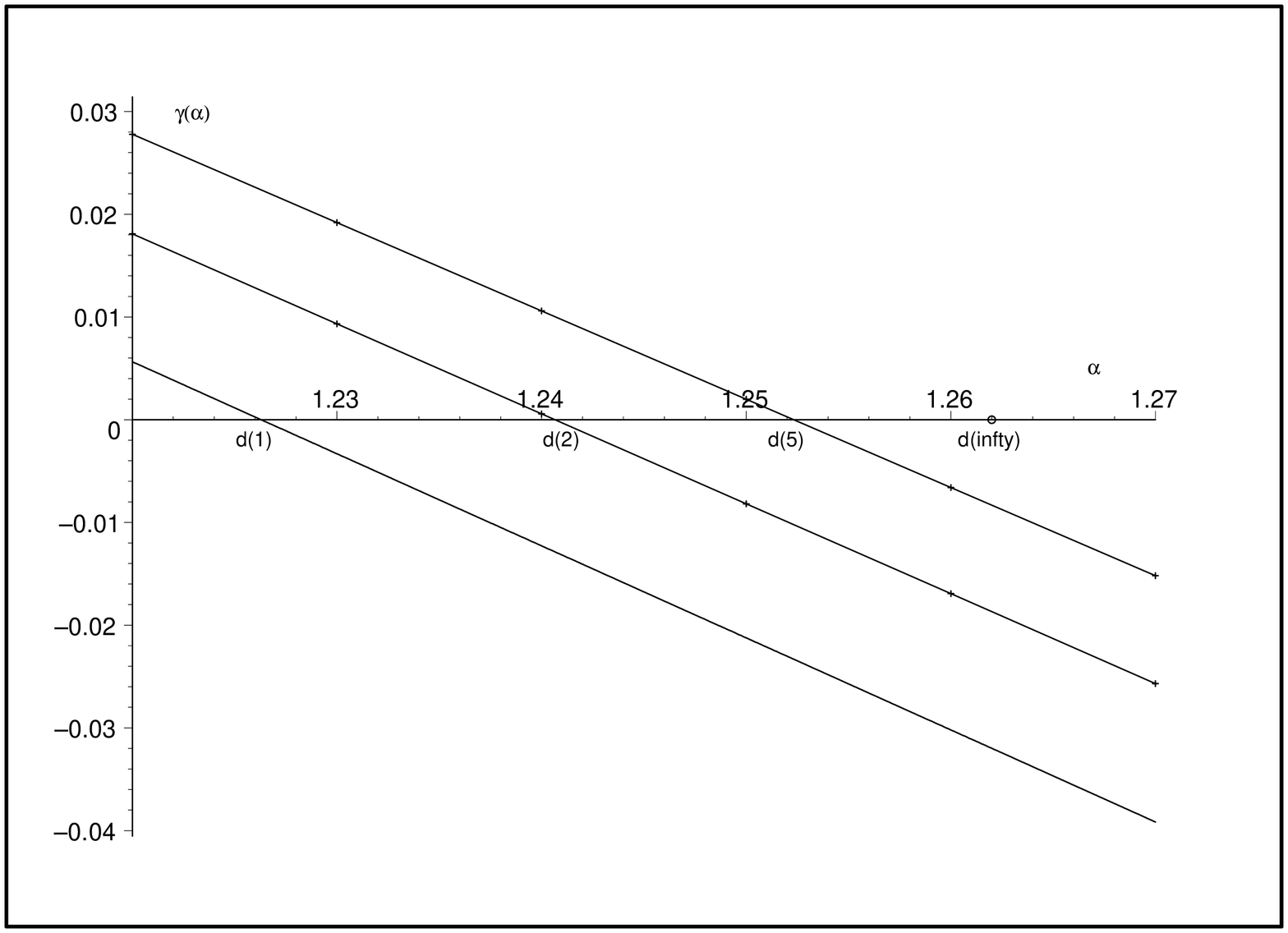} 
\end{center}
\begin{quote}
\textsf{\textbf{Figure  10.}  \textsf{Graphs of the ``pressure'' function  $ \gamma_V (\alpha ) $  for V = 1, 2, and 5 respectively, from left to right.}}
\end{quote}

 It was previously known that the dimension  of homogeneous random Sierpinski triangles 
 ($ V=1 $) is $ 2\log 3 /(\log 2 + \log 3) \approx 1.226 $, see Hambly (1992, 2000), and that the dimension of standard random Sierpinski triangles ($  V \to \infty$) is the solution $ d $ of $ \frac{1}{2} 3 \left(\frac{1}{2}\right)^d +  \frac{1}{2} 3 \left(\frac{1}{3}\right)^d=1 $, or approximately 1.262, see Falconer (1986), Graf (1987) and Mauldin and Williams  (1986).     
In particular, \eqref{pf1} is in agreement when computing $d(1) $. For  $ V > 1 $  we used Monte Carlo simulations to compute $ \gamma_V(\alpha) $   in the region of the interval $ [d(1), d(\infty)] $, with the computed values shown (Figure 10). These values have error at most .001 at the 95\% 
confidence level, and from this  one obtains the dimensions 
$ d(2) \approx 1.241 $, $ d(5) \approx 1.252 $ (Figure 10).  The computed graphs  for $ V>1 $  are concave up, although this does not show on the scale of Figure  10.

\section{Analysis  of  Dimension results for  $ V $-variable Fractals}

In order to motivate the following analysis, consider a smooth curve or smooth surface having  dimension 1 or 2 respectively.  It is possible to cover each ``efficiently'' (i.e.\ with little overlap) by sets of small diameter such that the sum of the diameters raised to the power 1 or 2 respectively is very close to the length, or to the area divided by $ \pi/4 $, respectively.  However, if any power $ \alpha > $  1 or 2 respectively is used then the limit of this sum, as the maximum diameter of the covering sets becomes arbitrarily small, is zero.  For any power $ \alpha  <$ 1 or 2 respectively the limit of the sum, as the maximum diameter becomes arbitrarily small, is infinity. 

In the case of the construction of $ 5 $-variable Sierpinski fractal triangles, we see that if one begins the construction process with 5 copies of a triangle $ T $  as indicated (Figure 3), 
then the contents of each buffer at later stages will consist of a large number of tiny triangles which approximate  and cover the ``ideal'' or limiting $ 5 $-variable Sierpinski fractal triangles.  (Note that what is obtained after $ k $ steps is only an actual $ 5 $-variable Sierpinski fractal triangle up to $ k $ levels of magnification.) One can check that the sum $ S^k_v(\alpha) $ of the diameters to the power $ \alpha  $ of the triangles in the $ v $-th buffer at level $ k $ is given by the sum of the entries in the $ v $-th row of the matrix $ M^1(\alpha)\cdot \ldots \cdot M^k(\alpha) $.

Using \eqref{pf} it is not too difficult to show that $  \lim_{k\to \infty} \frac{1}{k} \log S^k_v(\alpha) $ also exists for each $ v $ and equals $ \gamma_V(\alpha)$, independently of $ v $,  with probability one.  (The argument relies on the existence of ``necks'',  where a neck in the construction process occurs at some level if the same IFS and the same single fixed  input buffer is used for constructing the set  in each buffer at that level.) 
One can show (Barnsley, Hutchinson and Stenflo 2003b) that $ \gamma_V(\alpha ) $ is decreasing in $ \alpha $ and deduce that there is a unique $ d=d(V) $ such that $ \gamma_V(d)  = 0$ (Figure  10).  From this and \eqref{pf} it follows that
\begin{align*}
\alpha < d  \Longrightarrow  \gamma(\alpha) > 0 \Longrightarrow 
&\lim_{k\to \infty} S^k_v(\alpha) = \lim_{k\to \infty} \exp(k\, \gamma(\alpha))=\infty,\\
\alpha > d  \Longrightarrow  \gamma(\alpha) < 0 \Longrightarrow 
  &\lim_{k\to \infty} S^k_v(\alpha) = \lim_{k\to \infty} \exp(k\, \gamma(\alpha))=0.
\end{align*}
It is hence plausible from the previous discussion of curves and surfaces that the dimension of 
$ V $-variable Sierpinski triangles equals $d(V)  $ with probability one.  The motivation is that the covering by small triangles is very ``efficient''.

The justification that the dimension is at most $ d $ is in fact now straightforward from the definition of (Hausdorff) dimension of a set.  The rigorous argument that the dimension is at least $ d $, and hence exactly $ d $, is much more difficult, see Barnsley, Hutchinson and Stenflo (2003b).  It requires a careful analysis of the frequency of occurrence of necks and the construction of Gibbs type measures on $ V $-variable Sierpinski triangles, analogous to ideas in statistical mechanics. 

Similar results on dimension have been established much more generally,  see Barnsley, Hutchinson and Stenflo (2003b).  For example, 
suppose the functions in each IFS are similitudes, i.e.\  built up from translations, rotations, reflections in lines, and a single contraction around a fixed point by a fixed ratio $ r $ (both the point and the ratio $r  $  may depend on the  function in question).  
We also require that the IFSs involved satisfy the \emph{uniform open set condition}.  In the case  of $ 5 $-variable Sierpinski fractal triangles constructed from the IFSs $ F=(f_1,f_2,f_3) $ and $ G = (g_1,g_2,g_3) $ this means the following. There is an open set $ O $ (the interior of the triangle $ T $) such that  $ f_1(O) \subseteq O $, $ f_2(O) \subseteq O $, $ f_3(O) \subseteq O $, and 
$  f_1(O) $, $  f_2(O) $,  $  f_3(O) $ have no points in common, and such that analogous conditions apply to the maps $ g_1 $, $ g_2 $, $ g_3 $ with the \emph{same} set $ O $. 
Under these circumstances one constructs the matrices $ M^k(\alpha) $ and the pressure function $ \gamma_V(\alpha) $ as before and it follows   that the solution  $ d(V) $ of $ \gamma_V (d)=0 $ is the dimension of the corresponding $ V$-variable fractals with probability  one.

\section{Generalisations}
Many generalisations are possible.  

The number of functions  $ M $  in each of the IFSs  
$ F^n = \left( f^n_1,\dots, f^n_M;  p^n_1,\dots, p^n_M\right) $ may vary.  The maps  $f^n_m $   need only be mean contractive when their contraction ratios are averaged over  $ m $   and  $ n $  according to the probabilities  $ p^n_m $  and  $ P^n $  respectively.  Neither the number of functions in an IFS nor the number of IFSs need be finite; this is important for simulating various selfsimilar processes, including Brownian motion, see Hutchinson and R\"{u}schendorf (2000).  The maps  $f^n_m $ may be nonlinear, and many of the results and arguments, including those concerning dimension, will still be valid.  The set maps    $f^n_m $ need not be induced from point maps; this is technically useful in extending results to the case where $ M $   is not constant by artificially adding set maps  $f^n_m $   such that  $f^n_m(A) $   is always the empty set.  It could also be important in applications to modelling biological or physical phenomena where the objects under consideration are not naturally modelled as sets or measures.  Buffer sampling need not be uniform; buffers could be placed in a rectangular or other grid, and nearby buffers sampled with greater probability, in order to simulate various biological and physical phenomena.  

An IFS operates on $\mathbb{R}^2 $, or more generally on a compact metric space $(\mathbb{X},d)$, to produce a fractal set attractor; a weighted IFS produces a fractal measure attractor.   We have seen in this paper how a family of  IFSs operating on $(\mathbb{X},d)$, a probability distribution on this family of IFSs, and an integer parameter $ V $, can be used to generate a (super)IFS  operating in a natural way on  $( \mathbb{H}(\mathbb{X})^V, d_{\mathcal{H}})$, where 
$\mathbb{H}(\mathbb{X})$ is the space of compact subsets of $\mathbb{X}$ and $d_{\mathcal{H}}$ is the induced Hausdorff metric.  In the case of a family of weighted IFSs, the induced superIFS operates on  
$( \mathbb{P}(\mathbb{X})^V, d_{MK})$ where $\mathbb{P}(\mathbb{X})$ is the space of unit mass measures on $\mathbb{X}$ and $d_{MK}$ is the induced Monge Kantorovitch metric.  In either case there is a  superfractal set (consisting of $V$-variable sets or measures respectively) together with an associated superfractal measure (a probability distribution on the collection of $V$-variable sets or measures).  There is also a fast forward algorithm to generate this superfractal.

We can consider iterating this procedure.  Replace $(\mathbb{X},d)$ by 
$(\mathbb{H}(\mathbb{X})^V,d_{\mathcal{H}})$ or $( \mathbb{P}(\mathbb{X})^V, d_{MK})$, take a 
family of superIFSs, an associated probability distribution, and a new parameter $W$.
To speculate; if a superfractal may be thought of as a gallery of  a new class of fractal images, can one use some version of the iteration to produce a museum of galleries of yet another new class of fractal images?

\section{Conclusion}
There appear to be many potential applications, which include both the extension of modelling possibilities to allow a controlled degree of variability where deterministic or random fractals have been previously applied, and the rapid generation of accurately distributed examples of random fractals --- previously not possible except in very special cases.


\ack We thank the Australian Research Council for their support of this research, which was carried out at the Australian National University.

\section*{References}
\begin{harvard}

\item[] Barnsley M  F and  Demko S 1985  Iterated function systems and the global construction of fractals {\it Proc. Roy. Soc. Lond. Ser. A}  {\bf 399} 243--275 

\item[] Barnsley M  F, Hutchinson  J E and Stenflo  \"{O} 2003a A fractal value random iteration algorithm and fractal hierarchy {\it submitted}

\item[] Barnsley M  F, Hutchinson  J  E and Stenflo  \"{O} 2003b Dimension and approximation properties of $ V $-variable fractals {\it in preparation}

\item[] Cohen  J  E  1988  Generalised products of random matrices and operations research  {\it SIAM Review} {\bf 30} 69--86 

\item[] Falconer  K  J   1986 Random fractals {\it Math. Proc. Cambridge Philos. Soc.} {\bf 100} 
559--82 

\item[] Falconer K  J 1990  {\it Fractal Geometry. Mathematical Foundations and Applications}    (Chichester: John Wiley \& Sons)

\item[] Falconer K  J  1997 {\it Techniques in Fractal Geometry} (Chichester: John Wiley \& Sons) 

\item[] Feng  D J  and Lau  K S  2002 The pressure function for products of non-negative matrices  {\it Math. Res. Letters} {\bf 9} 1--16  

\item[] Furstenberg  H and Kesten  H 1960  Products of random matrices {\it Ann. Math. Statist. 31} 457--469 

\item[] Graf  S 1987   Statistically self-similar fractals {\it Probab. Theory Related Fields} 
{\bf 74} 357--392 

\item[] Hambly  B  M 2000  Heat kernels and spectral asymptotics for some random Sierpinski gaskets {\it Progr. Probab.}  {\bf 46} 239--267  
 
\item[] Hambly  B  M  1992  Brownian motion on a homogeneous random fractal {\it Probab. Theory Related Fields} {\bf 94} 1--38  

\item[] Hutchinson  J E 1981 Fractals and self-similarity {\it Indiana Univ. Math. J.}  {\bf  30}  713--749

\item[] Hutchinson J E  and R\"{u}schendorf  L 1998  Random fractal measures via the contraction method {\it Indiana Univ. Math. J.} {\bf  47}  471--487 
 
\item[] Hutchinson  J E and R\"{u}schendorf  L 2000  Random fractals and probability metrics {\it  Adv. in Appl. Probab.}  {\bf 32}  925--947 
 
\item[] Hutchinson  J E  and R\"{u}schendorf  L  2000  Selfsimilar fractals and selfsimilar random fractals {\it Progr. Probab.}  {\bf 46} 109--123 

\item[] Mauldin  R D  and Williams  S  C 1986  Random recursive constructions; asymptotic geometrical and topological Properties. {\it Trans. Amer. Math. Soc.} {\bf  295}
 325--346 
 
\item[] Stenflo  \"{O} 2001  Markov chains in random environments and random iterated function systems  {\it Trans. Amer. Math. Soc.} {\bf 353} 3547--62

\end{harvard}

\end{document}